\documentclass[3p]{elsarticle}

\usepackage{amssymb}
\usepackage{amsthm}

\usepackage{amsmath,amsthm,a4wide,mathrsfs,color,array}
\usepackage[utf8]{inputenc}
\usepackage[T1]{fontenc}
\usepackage{subcaption}
\usepackage{graphicx}
\usepackage{float}
\usepackage{multirow}

\newtheorem{theorem}{Theorem}
\newtheorem{proposition}{Proposition}
\newtheorem{remark}{Remark}
\newtheorem{assumption}{Assumption}
\newtheorem{lemma}{Lemma}

\newcommand{\nc}{\newcommand}
\nc{\dsp}{\displaystyle}
\nc{\mrm}{\mathrm}
\nc{\lbr}{\lbrack}
\nc{\rbr}{\rbrack}

\nc{\llbr}{\lbrack\!\lbrack}
\nc{\rrbr}{\rbrack\!\rbrack}

\nc{\mB}{\mrm{B}}
\nc{\mL}{\mrm{L}}
\nc{\bmL}{\boldsymbol{L}}
\nc{\mH}{\mrm{H}}
\nc{\mY}{\mrm{Y}}
\nc{\mV}{\mrm{V}}
\nc{\mK}{\mrm{K}}
\nc{\mS}{\mrm{S}}
\nc{\mT}{\mrm{T}}
\nc{\mbH}{\mathbb{H}}

\nc{\mM}{\mrm{M}}
\nc{\mD}{\mrm{D}}

\nc{\R}{\mathbb{R}}
\nc{\C}{\mathbb{C}}
\nc{\Z}{\mathbb{Z}}

\nc{\bx}{\boldsymbol{x}}
\nc{\by}{\boldsymbol{y}}
\nc{\bn}{\boldsymbol{n}}
\nc{\bu}{\boldsymbol{u}}
\nc{\bv}{\boldsymbol{v}}
\nc{\bvf}{\boldsymbol{f}}
\nc{\bp}{\boldsymbol{p}}
\nc{\bq}{\boldsymbol{q}}

\nc{\bL}{\boldsymbol{L}}

\nc{\ctru}{\mathsf{u}}
\nc{\ctrv}{\mathsf{v}}
\nc{\ctrf}{\mathsf{f}}

\nc{\bctru}{\boldsymbol{\mathsf{u}}}
\nc{\bctrv}{\boldsymbol{\mathsf{v}}}
\nc{\bctrf}{\boldsymbol{\mathsf{f}}}

\renewcommand{\div}{\mrm{div}}
\nc{\curl}{\mathbf{curl}}
\nc{\loc}{\mrm{loc}}
\nc{\inc}{\mrm{inc}}
\nc{\Id}{\mrm{Id}}
\nc{\MTFloc}{\mrm{MTF}_{\loc}}
\nc{\diag}{\mrm{diag}}

\nc{\IA}{\mathbb{A}}
\nc{\IB}{\mathbb{B}}
\nc{\II}{\mathbb{I}}
\nc{\IS}{\mathbb{S}}
\nc{\IK}{\mathbb{K}}

\nc{\imathj}{\mrm{i}}

\nc{\e}{\mrm{e}}

\nc{\bfP}{\mathbf{P}}
\nc{\bfQ}{\mathbf{Q}}
\nc{\bfW}{\mathbf{W}}
\nc{\bfV}{\mathbf{V}}
\nc{\bfU}{\mathbf{U}}
\nc{\bfE}{\mathbf{E}}

\newcommand{\abs}[1]{\left|#1\right|}
\renewcommand{\tan}{\textsc{t}}
\nc{\cur}{\textsc{r}}

\nc{\Hcurloc}[1]{\boldsymbol{H}(\curl,#1)}
\nc{\bfH}{\boldsymbol{H}}

\nc{\mG}{\mrm{G}}
\nc{\mA}{\mrm{A}}
\nc{\SL}{\mrm{SL}}
\nc{\DL}{\mrm{DL}}

\nc{\gr}{\textsc{g}}
\nc{\cu}{\textsc{c}}
\nc{\bfX}{\mathbf{X}}
\newcommand{\pp}{\parallel}

\nc{\ctrJ}{\mathfrak{J}}
\nc{\ctrH}{\mathfrak{H}}

\definecolor{dark_blue}{RGB}{46,87,144}
\newcommand\rev[1]{\textcolor{black}{#1}}

\usepackage[normalem]{ulem}

\journal{}

\begin{document}

\begin{frontmatter}



\title{Local Multiple Traces Formulation for Electromagnetics: Stability and Preconditioning for Smooth Geometries}
\tnotetext[mytitlenote]{This work received the support from ECOS-Conicyt under grant C15E07.
The first and second author would also like to acknowledge support from the French National
Research Agency (ANR) under grant ANR-15-CE23-0017-01. C.~Jerez-Hanckes thanks the support of Fondecyt Regular 1171491.
}
\author[1]{Alan Ayala}
\author[1]{Xavier Claeys}
\author[2]{Paul Escapil-Inchausp\'e}
\author[2]{Carlos Jerez-Hanckes\corref{cor1}}

\cortext[cor1]{Corresponding author \ead{carlos.jerez@uai.cl}}

\address[1]{Sorbonne Universit\'e, Universit\'e Paris-Diderot SPC, CNRS, Inria, Laboratoire Jacques-Louis Lions, \'equipe Alpines, Paris, France}
\address[2]{Faculty of Engineering and Sciences, Universidad Adolfo Ib\'a\~nez, Santiago,  Chile}

\begin{abstract}
We consider the time-harmonic electromagnetic transmission problem for the unit sphere. Appealing to a vector spherical harmonics analysis, we prove the first stability result  of the local multiple trace\rev{s} formulation (MTF) for electromagnetics, originally introduced by Hiptmair and Jerez-Hanckes [Adv.~Comp.~Math.~{\bf 37} (2012), 37-91] for the acoustic case, paving the way towards an extension to general piecewise homogeneous scatterers. Moreover, we investigate preconditioning techniques for the local MTF scheme and study the accumulation points of induced operators. In particular, we propose a novel second-order inverse approximation of the operator. Numerical experiments validate our claims and confirm the relevance of the preconditioning strategies.
\end{abstract}



\end{frontmatter}


\section{Introduction}\label{sec:Introduction}
Developing efficient computational methods for modeling
electromagnetic wave scattering by composite objects in unbounded
space remains a challenging problem raising many technical and
theoretical issues. Due to their rigorous account of radiation
conditions, boundary integral representations are among the preferred
choices despite the cumbersome electric and magnetic field integral
operators. However, when dealing with many subdomains, the problem can
become daunting computationally leading to high memory and CPU
requirements. Several boundary approaches have been proposed to tackle
electromagnetic wave transmission problems by homogeneous scatterers
in the frequency domain \cite{HJ12,CHJ13,CHJ15,CHH77,POM73,WTS77}. For
instance, the Poggio-Miller-Chang-Harrington-Wu and Tsai formulation
(PMCHWT), also referred to as Single Trace Formulation (STF), has been
shown to be useful for a number of applications and amenable to
significant improvements in terms of preconditioning---when using
iterative solvers---for the case of separated scatterers
\cite{CAF11,KBH18,KBH22}. For Laplace, Helmholtz and Maxwell's
equations, \emph{Multiple Traces Formulations} (MTFs)
\cite{HJ12,CHJ13} were introduced as a mean\rev{s} to solve
transmission problems by multiple connected scatterers and allowing to
use Calder\'on-based preconditioners. Though theoretical aspects for
Maxwell scattering are now fully available for \emph{global MTFs}
\cite{MR2996334,MR3069956} their implementation is extremely
cumbersome. Opposingly, local versions of the MTF
\cite{HJ12,HJL14,JPT15,HJA16,JPT17} are easily implemented and
parallelized. Although theoretical results for
acoustic and static versions are available, in the electromagnetic
case similar results remain elusive. Regarding Maxwell's equations,
preconditioning is crucial as STF and MTFs incorporate \emph{electric
field integral operators}. These commonly generate highly
ill-conditioned linear systems and lead to solver time being the
bottleneck of such schemes, or even to stagnation of iterative solvers
(refer e.g.~to \cite{EIJH19,KBH22}).

In the present contribution, we investigate various theoretical
aspects of the local MTF applied to Maxwell's equations. We will focus
on a geometrical configuration made up by two smooth subdomains and
one interface. A large part of this article actually assumes that the
interface is a sphere, which is still relevant as conclusions for
general smooth interfaces can be drawn from this particular case
arguing by compact perturbation, see e.g. \cite[Chap. 2]{SS11} or
\cite[Chap.5]{HW08}. In this canonical setting, one can directly use
separation of variables via vector spherical harmonics. One of our
main results is the derivation of a G\aa rding-type inequality for
the local MTF in the electromagnetic context (see Theorem
\ref{GardingIneq}).

We also study in detail the essential spectrum of the local MTF
operator and its preconditioned variants, not thoroughly studied
before to our knowledge. We exhibit surprisingly simple
formulas---\eqref{eq:eigInfty} and \rev{\eqref{eq:accPoints}}---for
these accumulation points.  Finally, we propose and analyze several
preconditioning techniques for the local MTF, looking for strategies
that (i) reduce as much as possible the number of accumulation points
in the spectrum of the preconditioned operator; and, (ii) lead to a
\emph{second-kind Fredholm operator} on smooth surfaces---a compact
perturbation of the identity operator.

The outline of this article is as follows. In Sections \ref{PbSetting}
and \ref{sec:TraceSpacesAndOperators} we set the problem under study
and introduce a necessary notation and definitions related to trace
spaces and potential theory. In Section \ref{sec:LocalMTF} we derive
the local MTF for Maxwell's equations. In Section
\ref{ProofOfUniqueness} we show that the kernel of this operator is
trivial.  In Section \ref{SeparationOfVariables} we apply separation of variables to the local MTF
operator and study the asymptotics of its spectrum. In Section
\ref{StabilityOfLocalMTF} we use separation of variables to establish
a G\aa rding-type inequality thus proving that the local MTF operator
is an isomorphism that, under conforming Galerkin discretization,
leads to quasi-optimal numerical methods. Section
\ref{sec:Preconditioning} describes several preconditioning strategies
along with numerical tests to discuss their performances. Finally,
concluding remarks are provided in Section \ref{sec:Conclusion}.

\section{Problem setting and functional spaces}
\label{PbSetting}
\begin{figure}[t]
  \centering
  \includegraphics[width=.3\linewidth]{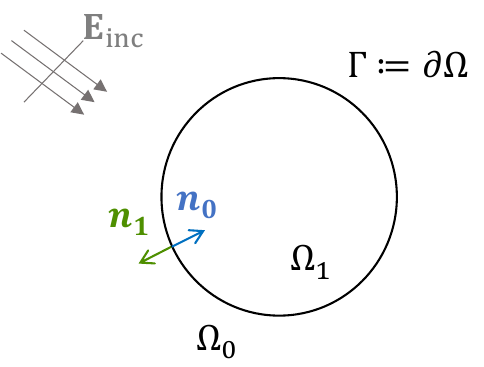}
\caption{Problem representation. Observe normal definitions.}
\label{fig:scatterer}
\end{figure}

We consider a partition of free space as $\R^{3} = \overline{\Omega}_{0}\cup \overline{\Omega}_{1}$ 
in\rev{to} two smooth open subdomains such that $\Omega_{0} = \R^{3}\setminus\overline{\Omega}_{1}$ (see Fig.~\ref{fig:scatterer}). \rev{For $j=0,1$, we set $\Gamma_j := \partial\Omega_{j}$ and $\Gamma := \Gamma_0 = \Gamma_1$. Let $\bn_{j}$} refer to the unit normal vector to $\Gamma$ directed toward the exterior of $\Omega_{j}$, so that we have $\bn_{0} = -\bn_{1}$. \rev{We define the \emph{skeleton} $\Sigma := \Gamma_0 \cup \Gamma_1$. }Let $\epsilon_{j}>0$ (resp.~$\mu_{j}>0$) refer to the electric permittivity (resp.~magnetic permeability) in the domain $\Omega_{j}$. 

We are interested in computing the scattering of an incident 
electromagnetic wave $(\mathbf{E}_{\inc},\mathbf{H}_{\inc})$ propagating in time-harmonic regime at pulsation $\omega>0$. To simplify matters, we require that $\curl(\mathbf{E}_{\inc})  -\imathj \omega \mu_{0}\mathbf{H}_{\inc} = 0$ and 
$\curl(\mathbf{H}_{\inc})  +\imathj \omega \epsilon_{0}\mathbf{E}_{\inc} = 0$ in $\R^{3}$: the incident field 
may, for example, be a plane wave. The equations for the total electromagnetic field $(\mathbf{E},\mathbf{H})$ under consideration read
\rev{\begin{subequations}
\label{MaxwellPb:Eq1}
\begin{align}
\curl(\mathbf{E})  -\imathj \omega \mu_{j}\mathbf{H} &= 0
&& \text{in}\ \Omega_j \text{ (for $j=1,2$)},\\
\curl(\mathbf{H})  +\imathj \omega \epsilon_{j}\mathbf{E}& = 0
&& \text{in}\ \Omega_j \text{ (for $j=1,2$)},\\
\abs{ \sqrt{\mu_{0}}(\mathbf{H} - \mathbf{H}_{\inc})\times \hat{\bx} - 
\sqrt{\epsilon_{0}}(\mathbf{E} - \mathbf{E}_{\inc}) }&= \mathcal{O}_{\vert\bx\vert\to\infty}(\vert \bx\vert^{-2}),\label{MaxwellPb:Eq1c}\\
\bn_{0}\times\mathbf{E}\vert_{\Gamma_{0}} + \bn_{1}\times\mathbf{E}\vert_{\Gamma_{1}}&  = 0 ,\label{MaxwellPb:Eq1d}\\
\bn_{0}\times\mathbf{H}\vert_{\Gamma_{0}} + \bn_{1}\times\mathbf{H}\vert_{\Gamma_{1}}&  = 0.\label{MaxwellPb:Eq1e}
\end{align}
\end{subequations}}
In addition, \rev{\eqref{MaxwellPb:Eq1c}} is referred to as the Silver-M\"uller radiation condition \cite{NED01,M03} wherein $\hat{\bx}:= \bx/\vert\bx\vert$. In  
\eqref{MaxwellPb:Eq1d}--\eqref{MaxwellPb:Eq1e} the notation "$\mathbf{E}\vert_{\Gamma_{j}}$" (resp. "$\mathbf{H}\vert_{\Gamma_{j}}$") should be 
understood as the trace taken at $\Gamma$ from the interior of $\Omega_{j}$---precise definitions provided in Section \ref{sec:TraceSpacesAndOperators}. \rev{For the sake of clarity, we represent the problem under consideration in Figure \ref{fig:scatterer}.}

One can reformulate \rev{\eqref{MaxwellPb:Eq1}} as a second order transmission boundary value problem, which is the basis of the Stratton-Chu potential theory. Specifically,
\rev{\begin{subequations}
\label{MaxwellPb:Eq3}
\begin{align}
\curl \ \curl(\bfE) - \kappa_{j}^{2}\bfE &= 0
&& \text{in}\ \Omega_j \text{ (for $j=1,2$)},\\
 \abs{\curl(\bfE-\bfE_{\inc})\times\hat{\bx} - \imathj\kappa_{0}(\bfE-\bfE_{\inc})} &= \mathcal{O}_{\vert\bx\vert\to\infty}(\vert \bx\vert^{-2}),\\
 \bn_{0}\times\mathbf{E}\vert_{\Gamma_{0}} + \bn_{1}\times\mathbf{E}\vert_{\Gamma_{1}}  &= 0,\label{MaxwellPb:Eq3c}\\
 \mu_{0}^{-1}\bn_{0}\times\curl(\mathbf{E})\vert_{\Gamma_{0}} + \mu_{1}^{-1}\bn_{1}\times\curl(\mathbf{E})\vert_{\Gamma_{1}} & = 0.\label{MaxwellPb:Eq3d}
\end{align}
\end{subequations}}

In the equations above we defined the effective wavenumber in each subdomain:
\begin{equation}\label{EffectiveWavenumber}
  \kappa_{j} := \omega\sqrt{\mu_{j}\epsilon_{j}}\quad j=0,1. 
\end{equation}
We study the solution of this problem by means of a boundary integral formulation. As mentioned before, \rev{several formulations are possible} but we focus here on local MTF.  As a complete stability analysis of the local MTF for the electromagnetic case is not presently available, we concentrate
on the following special case.

\begin{assumption} $\Omega_{1}$ is the unit ball and $\Gamma$ is the unit sphere.
\end{assumption}
This will allow for explicit calculus by means of separation of variables which will help investigate and clarify 
the structure of operators associated \rev{with} the local MTF.

\section{Trace spaces and operators}\label{sec:TraceSpacesAndOperators}
We refer to \cite{M03} for a detailed survey of vector functional spaces for Maxwell's equations. We introduce three interior trace operators taken from the \textit{interior} of $\Omega_{j}$ and are defined, for all $\bfU\in[\mathscr{C}^{\infty}(\R^{3})]^{3}$, the space of infinitely differentiable volume vector fields, as
\begin{equation}\label{TraceOperators}
\begin{aligned}
\gamma^{j}_{\tan}(\bfU) := \bfU\vert_{\Gamma}\times\bn_{j},\quad\gamma^{j}_{\cur}(\bfU) := \curl(\bfU)\vert_{\Gamma}\times\bn_{j}, \quad 
& \gamma^{j}(\bfU) := (\gamma^{j}_{\tan}(\bfU),\gamma^{j}_{\cur}(\bfU)).
\end{aligned}
\end{equation}
By density, one can show $\gamma_{\tan}^{j},\gamma_{\cur}^{j}:
\Hcurloc{\overline{\Omega}_{j}}\to \bfH^{-\frac{1}{2}}(\div,\Gamma_{j})$ and $\gamma^{j}:\Hcurloc{\overline{\Omega}_{j}}\to
\bfH^{-\frac{1}{2}}(\div,\Gamma_{j})^{2}$, where $\bfH^{-\frac{1}{2}}(\div,\Gamma)$ refers to the space of tangential traces of volume-based vector fields. The space $\bfH^{-\frac{1}{2}}(\div,\Gamma_{j})$ is put in duality with itself via the bilinear form:
$$(\bu,\bv)\mapsto \int_{\Gamma}\bn_{j}\cdot (\bu\times \bv) d\sigma.$$

The trace operators $\gamma^{j}_{\tan,c}$ (resp.~$\gamma^{j}_{\cur,c}$, $\gamma^{j}_{c}$) refer to exactly the
same operators as \eqref{TraceOperators} but with traces taken from the exterior along the same direction of 
$\bn_{j}$. Then, we shall define jump and averages traces as
\begin{equation*}
\begin{aligned}
& \{\gamma^{j}_{\star}\}(\bfU)     := \frac{1}{2}(\gamma^{j}_{\star}(\bfU) + \gamma^{j}_{\star,c}(\bfU))\rev{,}\\
& \lbr\gamma^{j}_{\star}\rbr(\bfU) :=  \gamma^{j}_{\star}(\bfU) - \gamma^{j}_{\star,c}(\bfU)\quad \text{for}\;\star = \tan,\cur,
\end{aligned}
\end{equation*}
and define $\{\gamma^{j}\}$ and $\lbr\gamma^{j}\rbr$ accordingly.

We also need to introduce
duality pairings for $\bfH^{-\frac{1}{2}}(\div,\Gamma)^{2} = \bfH^{-\frac{1}{2}}(\div,\Gamma)\times \bfH^{-\frac{1}{2}}(\div,\Gamma)$
defined by
\begin{equation}
\lbr (\bu,\bvf),(\bv,\bq)\rbr_{\Gamma_{j}}:=
\int_{\Gamma_{j}}\bn_{j}\cdot (\bu\times\bq + \bp\times \bvf)\, d\sigma.
\end{equation}
MTFs will be written in a so-called multiple traces space
and obtained as the Cartesian product of traces on the boundary of each
subdomain. In the present context, it takes the simple form:
\begin{equation}
\mbH(\Sigma):= \bfH^{-\frac{1}{2}}(\div,\Gamma)^{2}\times \bfH^{-\frac{1}{2}}(\div,\Gamma)^{2} = \bfH^{-\frac{1}{2}}(\div,\Gamma)^{4},
\end{equation}
\rev{with $\Sigma$ the skeleton.} This space will be equipped with a bilinear pairing $\llbr\cdot,\cdot\rrbr:\mbH(\Sigma)\times\mbH(\Sigma)\to \C$
defined as follows. For any tuples $\bctru = (\ctru_{0},\ctru_{1})$, $\bctrv = (\ctrv_{0},\ctrv_{1})$ in $\mbH(\Sigma)$,
we set
\begin{equation}
\llbr \bctru,\bctrv\rrbr := \lbr\ctru_{0},\ctrv_{0}\rbr_{\Gamma_{0}} + \lbr\ctru_{1},\ctrv_{1}\rbr_{\Gamma_{1}}.
\end{equation}
Note the identity $\llbr \bctru,\bctrv\rrbr = -\llbr \bctrv,\bctru\rrbr$ for any $\bctru,\bctrv\in \mbH(\Sigma)$.

\section{\rev{Local multiple traces operator}}\label{sec:LocalMTF}

As expected, we heavily rely on potential theory in the context of electromagnetics, i.e.~Stratton-Chu theory \cite{KH16}. In the sequel, let $\mathscr{G}_{\kappa}(\bx):=
\exp(\imathj\kappa\vert\bx\vert)/(4\pi\vert \bx\vert)$  refer to the outgoing
Green's kernel for the Helmholtz equation with wavenumber $\kappa>0$. 

Next, we define the boundary integral potentials: for $\ctru = (\bu,\bp)\in \bfH^{-\frac{1}{2}}(\div,\Gamma)^{2}$, we set
\begin{equation*}
\begin{aligned}
& \mG_{\kappa}(\ctru)(\bx):=\DL_{\kappa}(\bu)(\bx) + \SL_{\kappa}(\bp)(\bx)\\[5pt]
& \textrm{where}\quad \SL_{\kappa}(\bp)(\bx):=\kappa^{-2}
\int_{\Gamma}(\nabla\mathscr{G}_{\kappa})(\bx-\by)\div_{\Gamma}\bp(\by) d\sigma(\by) +
\int_{\Gamma}\mathscr{G}_{\kappa}(\bx - \by) \bp(\by)d\sigma(\by),\\
&  \textcolor{white}{\textrm{where}\quad}\DL_{\kappa}(\bu)(\bx):=
\curl\int_{\Gamma}\mathscr{G}_{\kappa}(\bx-\by) \bu(\by) d\sigma(\by).
\end{aligned}
\end{equation*}
The potential operator $\mG_{\kappa}$ maps continuously $\bfH^{-\frac{1}{2}}(\div,\Gamma)^{2}$ into $\bfH_{\loc}(\curl,\Omega_{0})$
and satisfies $(\curl\curl-\kappa_{0}^{2})\mG_{\kappa}(\ctru) = 0$ in $\Omega_{0}$ as well as Silver-M\"uller's radiation condition 
at infinity, regardless of $\ctru\in \bfH^{-\frac{1}{2}}(\div,\Gamma)^{2}$.  A similar result also holds in $\Omega_{1}$. The potential
operator plays a central role in the derivation of boundary integral equations as it can be used to represent solution to homogeneous Maxwell equations according to the Stratton-Chu integral representation theorem \cite[\rev{Theorem 5.49}]{KH16}. 

\begin{theorem}[{\rev{\cite[Theorem 5.49]{KH16}}}]\label{RepresentationFormula}
Let $\bfU\in \bfH_{\loc}(\curl,\Omega_{j})$ satisfy $\curl\curl (\bfU) - \kappa_{j}^{2}\bfU = 0$ in $\Omega_{j}$, $j=0,1$. For $j=0$
assume in addition that $\abs{\curl(\bfU)\times\hat{\bx} - \imathj\kappa_{0}\bfU}  = \mathcal{O}(\vert\bx\vert^{-2})$ 
for $\vert\bx\vert\to \infty$. Then we have $\mG_{\kappa}(\gamma^{j}(\bfU))(\bx) = 1_{\Omega_{j}}(\bx)\bfU(\bx)$ 
for all $\bx\in\R^{3}$.
\end{theorem}
On the other hand, the jumps of trace of the potential operator 
follow a simple and explicit expression.  
\begin{proposition}
\label{PotentialOperator} 
For any $\ctru\in \bfH^{-\frac{1}{2}}(\div,\Gamma)^{2}$ we have $\lbr\gamma^{j}\rbr\cdot\mG_{\kappa}(\ctru) = \ctru$.
\end{proposition}

In the forthcoming analysis, we shall make intensive use of the operator 
\rev{\begin{equation}\label{AOperator}
\mA_{\kappa}^{j}:= 2\{\gamma^{j}\}\cdot\mG_{\kappa}.
\end{equation}}
\rev{Standard choices of notation in the  literature dealing with Calder\'on projectors usually consider the same definition
  as above but without the factor 2. Our convention is motivated by simplifications stemming from this choice in our subsequent
  calculations.} It is clear that
$\{\gamma^{j}_{\tan}\}\cdot\DL_{\kappa} = \{\gamma^{j}_{\cur}\}\cdot\SL_{\kappa}$.
On the other hand, since the vector Helmholtz equation is satisfied by $\int_{\Gamma}\mathscr{G}_{\kappa}(\bx-\by) \bu(\by) d\sigma(\by)$, we find that that $\{\gamma^{j}_{\cur}\}\cdot\DL_{\kappa} = \kappa^{2}\{\gamma^{j}_{\tan}\}\cdot\SL_{\kappa}$.
As a consequence, the operator $\mA_{\kappa}^{j}$ can be represented in matrix form as
\begin{equation}\label{NotationOperator}
\mA_{\kappa}^{j}:=
\left\lbr\begin{array}{ll}
    \textcolor{white}{\kappa}\mrm{K}^{j}_{\kappa} & \kappa^{-1}\mrm{V}^{j}_{\kappa}\\
    \kappa\mrm{V}^{j}_{\kappa}                    & \textcolor{white}{\kappa^{-1}}\mrm{K}^{j}_{\kappa}
\end{array}\right\rbr
\quad\quad\text{where}\quad
\left\{\begin{array}{ll}
\mrm{V}^{j}_{\kappa} := \frac{2}{\kappa}\{\gamma^{j}_{\cur}\}\cdot\DL_{\kappa},\\[5pt]
\mrm{K}^{j}_{\kappa} := 2\{\gamma^{j}_{\tan}\}\cdot\DL_{\kappa}.
\end{array}\right.
\end{equation}
Observe that, for a given $\kappa$ we have $\mA_{\kappa}^{0} = -\mA_{\kappa}^{1}$ due to the change in the orientation of the normals $\bn_{0} = -\bn_{1}$. The operators (\ref{NotationOperator}) can be used to caracterise solutions of Maxwell's equations in a given subdomain. 
\begin{proposition}\label{CalderonProj}
The operator $\gamma^{j}\mG_{\kappa} = \frac{1}{2}(\Id+\mA_{\kappa}^{j})$ is a continuous projector as a mapping
from $\bfH^{-\frac{1}{2}}(\div,\Gamma)^{2}$ into $\bfH^{-\frac{1}{2}}(\div,\Gamma)^{2}$. Its range is the space of traces $\gamma^{j}(\bfU)$ where 
$\bfU\in\bfH_{\loc}(\curl,\Omega_{j})$ satisfies $\curl\curl(\bfU)-\kappa^{2}\bfU = 0\;\text{in}\;\Omega_{j}$, as well as
$\abs{\curl(\bfU)\times\hat{\bx} - \imathj\kappa_{0}\bfU}  = \mathcal{O}(\vert\bx\vert^{-2})$ for $\vert\bx\vert\to \infty$ if $j=0$.
\end{proposition}

An immediate consequence of the above proposition, \rev{combined with the notational convention \eqref{AOperator}
(in particular the mulitplicative factor 2 in there),} is that $(\mA_{\kappa}^{j})^{2} = \Id$, known as Calder\'on's identity.
As the incident field is a solution to Maxwell's equations with wavenumber $\kappa_{0}$ on $\R^{3}$---which includes $\Omega_{1}$---, then
so that $(\mA^{1}_{\kappa_{0}} - \Id)\gamma^{1}(\bfE_{\inc}) = 0$ according to the proposition above. Since on the other hand, it holds that $\mA^{1}_{\kappa_{0}} = -\mA^{0}_{\kappa_{0}}$ 
and $\gamma^{0}(\bfE_{\inc}) = -\gamma^{1}(\bfE_{\inc})$ (continuity of $\bfE_{\inc}$ across interfaces), we conclude that 
$\mA^{0}_{\kappa_{0}}\gamma^{1}(\bfE_{\inc}) =  - \gamma^{1}(\bfE_{\inc})$. Using Proposition \ref{CalderonProj}, we also see that  equation (\ref{MaxwellPb:Eq3}) can be reformulated as  $(\mA_{\kappa_{1}}^{1} - \Id)\gamma^{1}(\bfE) = 0$ on the one hand, and 
$(\mA_{\kappa_{0}}^{0} - \Id)(\gamma^{0}(\bfE) - \gamma^{0}(\bfE_{\inc})) = 0$ on the other hand. The latter is equivalent to 
$(\mA_{\kappa_{0}}^{0} - \Id)\gamma^{0}(\bfE) =  - 2\gamma^{0}(\bfE_{\inc})$.

Next, we need to reformulate the transmission conditions \rev{\eqref{MaxwellPb:Eq3c}--\eqref{MaxwellPb:Eq3d}}. Since these conditions are weighted 
by the permeability coefficients $\mu_{j}$, we introduce scaling operators:
$$\tau_{\alpha}:\bfH^{-\frac{1}{2}}(\div,\Gamma)^{2}\to \bfH^{-\frac{1}{2}}(\div,\Gamma)^{2},$$ 
defined by $\tau_{\alpha}(\bv,\bq):=(\bv,\alpha\bq)$. By the definition of the effective wavenumber in \eqref{EffectiveWavenumber}, we see that $\omega\mu/\kappa = \sqrt{\mu/\epsilon}$ and we can define the scaled operators:
\begin{equation}
\mA_{\kappa,\mu}^{j}:=\tau_{\omega\mu}^{-1}\cdot\mA_{\kappa}^{j}\cdot\tau_{\omega\mu} =
\left\lbr\begin{array}{ll}
    \textcolor{white}{\sqrt{\dsp{\epsilon/\mu}}}\;\mrm{K}^{j}_{\kappa}   & \sqrt{\dsp{\mu/\epsilon}}\;\mrm{V}^{j}_{\kappa}\\
    \sqrt{\dsp{\epsilon/\mu}}\;\mrm{V}^{j}_{\kappa}                      & \textcolor{white}{\sqrt{\dsp{\mu/\epsilon}}}\;\,\mrm{K}^{j}_{\kappa}
\end{array}\right\rbr.
\end{equation}
With this definition, we have $(\mA_{\kappa,\mu}^{j})^{2} = \Id$. 

The transmission conditions then are rewritten as
\begin{equation}
\tau_{\omega\mu_{0}}^{-1}\gamma^{0}(\bfE) +  \tau_{\omega\mu_{1}}^{-1}\gamma^{1}(\bfE) = 0.
\end{equation}
For the sake of conciseness, we will thus choose $\ctru_{j} = \tau_{\omega\mu_{j}}^{-1}\gamma^{j}(\bfE)$
as unknowns of our problem. As a consequence, \rev{\eqref{MaxwellPb:Eq3} can be cast as
\begin{subequations}
\label{IntegralSystem}
\begin{align}
(\mA_{\kappa_{0},\mu_{0}}^{0} - \Id)\ctru_{0}& =  - 2\tau_{\omega\mu_{0}}^{-1}\gamma^{0}(\bfE_{\inc}),\\
(\mA_{\kappa_{1},\mu_{1}}^{1} - \Id)\ctru_{1} &= 0,\\
\ctru_{0} +\ctru_{1} &= 0.
\end{align}
\end{subequations}
Now, let us rewrite (\ref{IntegralSystem})
in a matrix form. We first introduce the continuous map $\IA_{(\kappa,\mu)}:\mbH(\Sigma)\to \mbH(\Sigma)$ as a block diagonal operator $\IA_{(\kappa,\mu)}(\bctru):= (\mA_{\kappa_{0},\mu_{0}}^{0}(\ctru_{0}), \mA_{\kappa_{1},\mu_{1}}^{1}(\ctru_{1}))$ for any $\bctru = (\ctru_{0},\ctru_{1})\in\mbH(\Sigma)$, with subscript $(\kappa,\mu)$ representing the dependence on $(\kappa_0,\kappa_1, \mu_0,\mu_1)$.
The first two rows of (\ref{IntegralSystem}) can be rewritten as 
\begin{equation}\label{MatrixCalderon}
(\IA_{(\kappa,\mu)}-\Id)\bctru = \bctrf
\end{equation} 
where  $\bctru = (\ctru_{0},\ctru_{1})$ and $\bctrf = (- 2\tau_{\omega\mu_{0}}^{-1}\gamma^{0}(\bfE_{\inc}), 0)$.}
To enforce transmission conditions, we also need to consider an operator $\Pi:\mbH(\Sigma)\to \mbH(\Sigma)$
whose action consists in swapping traces from both sides of the interface. It is defined by
$\Pi(\ctru_{0},\ctru_{1}) := (\ctru_{1},\ctru_{0})$ for both $\ctru_{0}$ and $\ctru_{1}$ in $\bfH^{-\frac{1}{2}}(\div,\Gamma)^2$, so that
transmission conditions simply rewrite $\bctru = -\Pi(\bctru)$. Plugging the transmission operator into (\ref{MatrixCalderon}) 
leads to the local MTF of \eqref{MaxwellPb:Eq3}:
\rev{
\begin{equation}
\text{Find}\;\bctru\in \mbH(\Gamma)\;\text{such that} \; \mrm{MTF}_{\loc}(\bctru) = \bctrf
\end{equation}
where 
\begin{equation}\label{MTFLoc}
\mrm{MTF}_{\loc}:= \IA_{(\kappa,\mu)} +\Pi = 
\left\lbr\begin{array}{cc}
\mA_{\kappa_{0},\mu_{0}}^{0} & \Id\\[5pt]
\Id & \mA_{\kappa_{1},\mu_{1}}^{1}
\end{array}\right\rbr.
\end{equation}}
\section{Injectivity of local MTF for one subdomain}\label{ProofOfUniqueness}
We now prove the injectivity of the operator $\mrm{MTF}_{\loc}$ introduced above.
Assume that $\bctru = (\ctru_{0},\ctru_{1})\in \mbH(\Sigma)$ satisfies $\bctru\in \mrm{Ker}(\mrm{MTF}_{\loc})$ \rev{i.e.}
\begin{equation}
\label{ZeroSystem}
\left\lbr\begin{array}{cc}
\mA_{\kappa_{0},\mu_{0}}^{0} & \Id\\[5pt]
\Id & \mA_{\kappa_{1},\mu_{1}}^{1}
\end{array}\right\rbr\left\lbr\begin{array}{cc}
\ctru_{0}\\[5pt]
\ctru_{1}
\end{array}\right\rbr =0.
\end{equation}
In accordance with Theorem \ref{RepresentationFormula}, we define the (radiating) solution
$\bfU(\bx) := \mG_{\kappa_j}(\tau_{\omega\mu_{j}}(\ctru_{j}))(\bx)$ for $\bx\in \Omega_{j}, j=0,1$.
Taking interior traces, scaling both formulas, and using that $\ctru$ solves (\ref{ZeroSystem}) yields:
\begin{align*}
\tau_{\omega\mu_{0}}^{-1} \gamma^0(\bfU)&= (\Id+\mA^{0}_{\kappa_0,\mu_0})\ctru_0=\ctru_0 - \ctru_1,\\
\tau_{\omega\mu_{1}}^{-1} \gamma^1(\bfU)&= (\Id+\mA^{1}_{\kappa_1,\mu_1})\ctru_1=\ctru_1 - \ctru_0,
\end{align*}
hence the trace jump $\tau_{\omega\mu_{0}}^{-1}\gamma^{0}(\bfU) +  \tau_{\omega\mu_{1}}^{-1}\gamma^{1}(\bfU) = 0$,
leading to the conclusion that $\bfU$ is a Maxwell solution over the whole $\mathbb{R}^3$. By uniqueness of
the Maxwell radiating solution \cite{NED01,M03}, it holds that $\bfU\equiv0$, and so 
\begin{equation*}
  \begin{aligned}
    & \ctru_0 - \ctru_1 = 0,\\
    & (\Id+\mA^{0}_{\kappa_0,\mu_0})\ctru_0 = 0,\\
    & (\Id+\mA^{1}_{\kappa_1,\mu_1})\ctru_1 = 0.
  \end{aligned}
\end{equation*}
Thus, $\ctru_j$ is Cauchy data in $\Omega_j^c := \R^{3}\setminus\overline{\Omega}_{j}$ for wavenumber $\kappa_j$.
Set $\bfU^c(\bx) = (-1)^{j}\mG_{\kappa_j}(\tau_{\omega\mu_{j}}(\ctru_{j}))(\bx)$ for $\bx\in \Omega_{j}^{c}, j=0,1$,
and repeating the same arguments as above yields 
\begin{align*}
\tau_{\omega\mu_{0}}^{-1} \gamma_{c}^0(\bfU^c)&= (-\Id+\mA^{0}_{\kappa_0,\mu_0})\ctru_0= -(\ctru_0 + \ctru_1),\\
\tau_{\omega\mu_{1}}^{-1} \gamma_{c}^1(\bfU^c)&=(+\Id-\mA^{1}_{\kappa_1,\mu_1})\ctru_1= +(\ctru_1 + \ctru_0),
\end{align*}
giving $\tau_{\omega\mu_{0}}^{-1}\gamma_{c}^{0}(\bfU^{c}) +  \tau_{\omega\mu_{1}}^{-1}\gamma_{c}^{1}(\bfU^{c}) = 0$.
We see that $\bfU^{c}$ is solution to a one-subdomain transmission problem with homogeneous source term.
Such a problem admits zero as unique solution which yields $\bfU^{c} \equiv 0$. Finally,
$\tau_{\omega\mu_{j}}(\ctru_{j}) = \lbr\gamma^{j}\rbr\cdot \mG_{\kappa_j}(\tau_{\omega\mu_{j}}(\ctru_{j})) =
\gamma^{j}(\bfU)-(-1)^{j}\gamma^{j}_{c}(\bfU^{c}) = 0$ leading to $\ctru_{0} = \ctru_{1} = 0$. 
We have established the following \rev{result}.

\begin{proposition}
  Under the assumptions of Section \ref{PbSetting} we have $\mrm{Ker}(\mrm{MTF}_{\loc}) = \{{\bf 0}\}$.
\end{proposition}
  
\section{Spectral analysis of the local MTF operator for Maxwell equations}\label{SeparationOfVariables}
We are interested in deriving an explicit expression of operator (\ref{MTFLoc}) and analyzing the 
eigenvalues of the preconditioned formulations. As the present geometrical setting is spherically symmetric, this can be obtained by means of separation of variables based on spherical harmonics.

\subsection{Tangential spherical harmonics}
 Any tangential vector field
\rev{$$
\bu\in\bmL^{2}_{\tan}(\Gamma):=\left\{\bv:\Gamma\to \C \; \Big{|} \; \forall \bx \in \Gamma , \; \bv(\bx)\cdot\bx = 0 \quad \text{and} \quad  \Vert \bv\Vert_{\bmL^{2}_{\tan}(\Gamma)}^{2} := \int_{\Gamma}\vert\bv\vert^{2} d\sigma<+\infty \right\}
$$}
can be decomposed as
\cite{MR3268844,MR1822275}
\begin{equation}\label{HilbertExpansion}
\begin{aligned}
&\bu(\bx) = \sum_{n=0}^{+\infty}\sum_{\vert m\vert\leq n}
u_{n,m}^{\pp}\bfX_{n,m}^{\pp}(\bx)  + u_{n,m}^{\times}\bfX_{n,m}^{\times}(\bx)\\
&\text{with}
\quad\bfX_{n,m}^{\pp}:=\frac{1}{\sqrt{n(n+1)}}\nabla_{\Gamma}\mY_{n}^{m},
\quad\bfX_{n,m}^{\times}:=\bn_{1}\times\bfX_{n,m}^{\pp},
\end{aligned}
\end{equation}
where $\nabla_{\Gamma}$ is the surface gradient. Denoting $(\theta,\varphi)\in\lbr0,\pi\rbr\times \lbr0,2\pi\rbr$
the spherical coordinates on $\Gamma$, spherical harmonics (see e.g \cite[\S 14.30]{MR2723248}) are defined by
\begin{equation}
\mY_{n}^{m}(\theta,\varphi):=\sqrt{\frac{2n+1}{4\pi}}\sqrt{\frac{(n-\vert m\vert)!}{(n+\vert m\vert)!}}
\;\mrm{P}_{n}^{\vert m\vert}(\cos\theta)\exp(\imathj m\varphi).
\end{equation}
In the definition above, the functions $\mrm{P}_{n}^{m}(t), m\geq 0, t\in \lbr 0,1\rbr$ refer
to the associated Legendre functions, see e.g. \cite[\S7.12]{MR0350075}. The tangent fields
$\bfX_{n,m}^{\pp}, \bfX_{n,m}^{\times}$ form an orthonormal Hilbert basis
of $\mL^{2}_{\tan}(\Gamma)$. Let us denote $\bfX_{n,m}(\bx) := \lbr \bfX_{n,m}^{\pp}(\bx),
\bfX_{n,m}^{\times}(\bx)\rbr$ which maps a pair of scalar coefficients to a tangent vector
field over $\Gamma$ (\ref{HilbertExpansion}) can then be rewritten in the more compact form
\begin{equation}
\bu(\bx) = \sum_{n=0}^{+\infty}\sum_{\vert m\vert\leq n}\bfX_{n,m}(\bx)\cdot u_{n,m}
\end{equation}
for a collection of coordinate vectors $u_{n,m} = \lbr u_{n,m}^{\pp},u_{n,m}^{\times}\rbr^{\top}\in\C^{2}$.

\subsection{Local MTF operator over a sphere: separation of variables}

The operators coming into play in the expression of the local multi-trace operator (\ref{MTFLoc})
are actually (block) diagonalized by this basis. Define $\mathfrak{J}_{n}(t):= \sqrt{\pi t/2}J_{n+1/2}(t)$ where $J_{\nu}(t)$ are Bessel functions of the first kind of order $\nu$ (see \cite[\S10.2]{MR2723248})
and $\mathfrak{H}_{n}(t):=\sqrt{\pi t/2}H_{n+1/2}^{(1)}(t)$ where $H_{\nu}^{(1)}(t)$ are Hankel functions of the first kind of order $\nu$ (see \cite[\S10.2 \& \S10.4]{MR2723248}). Then, according to Lemma 1 in \cite{MR3268844} and
using notations (\ref{NotationOperator}), we have
\begin{equation}\label{MatrixForm1}
\begin{aligned}
  & \rev{\mV^{0}_{\kappa}\cdot\bfX_{n,m} = \mV^{0}_{\kappa}\cdot\lbr\bfX_{n,m}^{\pp}, \bfX_{n,m}^{\times}\rbr}\\
  & \textcolor{white}{\mV^{0}_{\kappa}\cdot\bfX_{n,m}} \rev{= \lbr\mV^{0}_{\kappa}\cdot\bfX_{n,m}^{\pp}, \mV^{0}_{\kappa}\cdot\bfX_{n,m}^{\times}\rbr
  = \bfX_{n,m}\cdot\widehat{\mV}^{0}_{\kappa}\lbr n\rbr}
  \quad \text{where}\\[5pt]
& \rev{\widehat{\mV}}^{0}_{\kappa}\lbr n\rbr :=
\left\lbr\begin{array}{cc}
0 & +2\imathj\,\mathfrak{J}_{n}(\kappa)\mathfrak{H}_{n}(\kappa)\\
-2\imathj\;\mathfrak{J}_{n}'(\kappa)\mathfrak{H}_{n}'(\kappa) & 0
\end{array}\right\rbr \;\;\in\C^{2\times 2}.
\end{aligned}
\end{equation}
Since $\mV^{1}_{\kappa} = -\mV^{0}_{\kappa}$ we have
$\mV^{1}_{\kappa}\cdot\bfX_{n,m}(\bx) = \bfX_{n,m}(\bx)\cdot\rev{\widehat{\mV}}^{1}_{\kappa}\lbr n\rbr$
by setting $\rev{\widehat{\mV}}^{1}_{\kappa}\lbr n\rbr := -\rev{\widehat{\mV}}^{0}_{\kappa}\lbr n\rbr$.
According to Lemma 1 in \cite{MR3268844}, we also have the explicit expression:
\begin{equation}\label{MatrixForm2}
\begin{aligned}
& \mK^{0}_{\kappa}\cdot\bfX_{n,m} = \bfX_{n,m}\cdot\rev{\widehat{\mK}}^{0}_{\kappa}\lbr n\rbr\quad \text{where}\\[5pt]
& \rev{\widehat{\mK}}^{0}_{\kappa}\lbr n\rbr := \imathj(\,\ctrJ_{n}(\kappa)\ctrH_{n}'(\kappa) + \ctrJ_{n}'(\kappa)\ctrH_{n}(\kappa)\,)
\left\lbr\begin{array}{cc}
-1 & 0\\
0 & +1
\end{array}\right\rbr \;\;\in\C^{2\times 2}.
\end{aligned}
\end{equation}
Here again, defining $\rev{\widehat{\mK}}^{1}_{\kappa}\lbr n\rbr = -\rev{\widehat{\mK}}^{0}_{\kappa}\lbr n\rbr$ we obtain $\mK^{1}_{\kappa}\cdot\bfX_{n,m}(\bx) = \bfX_{n,m}(\bx)\cdot\rev{\widehat{\mK}}^{1}_{\kappa}\lbr n\rbr$. From
(\ref{MatrixForm1}) and (\ref{MatrixForm2}) we deduce an explicit expression for the
operators $\mA_{\kappa,\mu}^{j}$. First of all define the function $\bfX_{n,m}^{\#\,2}$ by the expression
\begin{equation*}
\rev{  \bfX_{n,m}^{\#\,2}(\bx):=
  \left\lbr \begin{array}{cc}
    \bfX_{n,m}(\bx) & 0\\
    0 & \bfX_{n,m}(\bx)
  \end{array}\right\rbr =
  \left\lbr \begin{array}{cccc}
    \bfX_{n,m}^{\pp} & \bfX_{n,m}^{\times} & 0 & 0\\
    0 & 0 & \bfX_{n,m}^{\pp} & \bfX_{n,m}^{\times}
  \end{array}\right\rbr}
\end{equation*}
\rev{which we can also simply denote $\bfX_{n,m}^{\#\,2} = \mrm{diag}(\bfX_{n,m},\bfX_{n,m})$.
Then $\bfX_{n,m}^{\#\,2}$ should be understood as a linear operator that maps an element
of $\C^4$ to a pair of tangent vector fields over $\Gamma$. }
Any element $\ctru = (\bu,\bp)\in \bfH^{-\frac{1}{2}}(\div,\Gamma)^{2}$ decomposes as
$\ctru(\bx) = \sum_{n,m}\bfX_{n,m}^{\#\,2}(\bx)\cdot \ctru_{n,m}$
where $\ctru_{n,m}\in\C^{4}$ are coordinate vectors that do not depend on $\bx$.
In this basis, the operator $\mA_{\kappa,\mu}^{j}$ admits the following matrix form
\begin{equation}
\begin{aligned}
& \mA_{\kappa,\mu}^{j}\cdot \bfX_{n,m}^{\#\,2}(\bx) =
\bfX_{n,m}^{\#\,2}(\bx)\cdot \rev{\widehat{\mA}}_{\kappa,\mu}^{j}\lbr n\rbr\quad\text{where}\\[10pt]
& \rev{\widehat{\mA}}_{\kappa,\mu}^{j}\lbr n\rbr:=
\left\lbr\begin{array}{cc}
\textcolor{white}{\sqrt{\epsilon/\mu}}\;\widehat{\mK}^{j}_{\kappa}\lbr n\rbr & \sqrt{\mu/\epsilon}\;\widehat{\mV}^{j}_{\kappa}\lbr n\rbr\\[5pt]
\sqrt{\epsilon/\mu}\;\widehat{\mV}^{j}_{\kappa}\lbr n\rbr & \textcolor{white}{\sqrt{\mu/\epsilon}}\;\widehat{\mK}^{j}_{\kappa}\lbr n\rbr
\end{array}\right\rbr \;\;\in\C^{4\times 4}.
\end{aligned}
\end{equation}
\rev{We can reiterate the notational process used above,
and introduce the field $\bfX_{n,m}^{\#\,4}:=\mrm{diag}(\bfX_{n,m}^{\#\,2},\bfX_{n,m}^{\#\,2})$
which also writes in matrix form
\begin{equation*}
  \bfX_{n,m}^{\#\,4}(\bx) =
  \left\lbr\begin{array}{ll}
  \bfX_{n,m}^{\#\,2}(\bx) & 0 \\ 0 & \bfX_{n,m}^{\#\,2}(\bx)
  \end{array}\right\rbr.
\end{equation*}}
Then, any element $\bctru = (\ctru^{0},\ctru^{1})\in \bfH^{-\frac{1}{2}}(\div,\Gamma)^{2}\times \bfH^{-\frac{1}{2}}(\div,\Gamma)^{2}$
can be decomposed as $\bctru(\bx) = \sum_{n,m}\bfX_{n,m}^{\#\,4}(\bx)\cdot\ctru_{n,m}$ where $\ctru_{n,m}\in\C^{8}$ are coordinate
vectors that do not depend on $\bx$. Then the multi-trace operator (\ref{MTFLoc}) is reduced to matrix form in this basis
\begin{equation}
\begin{aligned}
& \mrm{MTF}_{\loc}\cdot \bfX_{n,m}^{\#\,4}(\bx)=\bfX_{n,m}^{\#\,4}(\bx)\cdot \rev{\widehat{\mrm{MTF}}}_{\loc}\lbr n\rbr\quad\text{where}\\[10pt]
& \rev{\widehat{\mrm{MTF}}}_{\loc}\lbr n\rbr:=
\left\lbr\begin{array}{cc}
\widehat{\mA}_{\kappa_{0},\mu_{0}}^{0}\lbr n\rbr & \Id\\[5pt]
\Id & \widehat{\mA}_{\kappa_{1},\mu_{1}}^{1}\lbr n\rbr
\end{array}\right\rbr \;\;\in\C^{8\times 8}
\end{aligned}
\end{equation}

\subsection{Accumulation points}\label{Asymptotics}
We can now study in more detail the symbol of the boundary integral operators introduced
in the previous section. To be more precise, we examine their behaviour for $n\to +\infty$.
First of all, from the series expansion of spherical Bessel functions given by \cite[\S10.53]{MR2723248},
we deduce that, for any fixed $t>0$, it holds that
\begin{equation}
\begin{aligned}
& \ctrJ_{n}(t) = t^{n+1}\frac{n!\,2^{n}}{(2n+1)!}\left\{1-\frac{t^{2}}{4n}+\mathcal{O}\left(\frac{1}{n^{2}}\right)\;\right\}\\
& \ctrH_{n}(t) = -\imathj t^{-n}\frac{(2n)!}{n!\,2^{n}}\left\{1+\frac{t^{2}}{4n}+\mathcal{O}\left(\frac{1}{n^{2}}\right)\;\right\}
\end{aligned}
\end{equation}
Since Bessel functions are expressed in terms of convergent series of analytic functions,
we can derive the above asymptotics. This leads to the following behaviours for
the derivatives,
\begin{equation}
\begin{aligned}
& \ctrJ_{n}'(t) = t^{n}\frac{n!\,2^{n}}{(2n+1)!}\left\{n+1-\frac{t^{2}}{4}+\mathcal{O}\left(\frac{1}{n}\right)\;\right\}\\
& \ctrH_{n}'(t) = \imathj t^{-(n+1)}\frac{(2n)!}{n!\,2^{n}}\left\{n+\frac{t^{2}}{4}+\mathcal{O}\left(\frac{1}{n}\right)\;\right\}
\end{aligned}
\end{equation}
One can combine these asymptotics to obtain the predominant behaviour of the functions coming into
play in the boundary integral operators expressions of the previous section. Specifically, we find
\begin{equation}
\begin{aligned}
& -2\imathj \ctrJ_{n}(t) \ctrH_{n}(t)  \mathop{\sim}_{n\to \infty} -t/n, \\
& +2\imathj \ctrJ_{n}'(t)\ctrH_{n}'(t) \mathop{\sim}_{n\to \infty} -n/t, \\
& \imathj (\;\ctrJ_{n}'(t) \ctrH_{n}(t)+ \ctrJ_{n}(t) \ctrH_{n}'(t)\;) \mathop{\sim}_{n\to \infty} 1/(2n).
\end{aligned}
\end{equation}
Define $\mT_{n}\in\mathbb{C}^{2\times 2}$ by $\mT_{n}(u_{1},u_{2}) := (u_{1}, u_{2}/n)$.
From this, we conclude that, as $n\to +\infty$, we have $\rev{\widehat{\mK}}^{0}_{\kappa}\lbr n\rbr\sim(2n)^{-1}\mK_{\kappa}^{0,\infty}$
and $\rev{\widehat{\mV}}^{0}_{\kappa}\lbr n\rbr\sim \mT_{n}^{-1}\cdot\mV_{\kappa}^{0,\infty}\cdot\mT_{n}$ where
$\mV_{\kappa}^{0,\infty},\mK_{\kappa}^{0,\infty}\in \mathbb{C}^{2\times 2}$ are constant matrices independent of $n$ given by
\begin{equation*}
\mV_{\kappa}^{0,\infty} :=
\left\lbr\begin{array}{cc}
0 & \kappa\\
1/\kappa & 0
\end{array}\right\rbr\quad\text{and}\quad
\mK^{0,\infty}_{\kappa}:=
\left\lbr\begin{array}{cc}
-1 &  0\\
0  & +1
\end{array}\right\rbr.
\end{equation*}
Next, define $\mT_{n}^{\#2}\in \mathbb{C}^{4\times 4}$ by $\mT_{n}^{\#2}(\bu_{1},\bu_{2}) = (\mT_{n}(\bu_{1}),\mT_{n}(\bu_{2}))$
for any pair $\bu_{1},\bu_{2}\in\mathbb{C}^{2}$  i.e. $\mT_{n} = \mrm{diag}(\mT_{n},\mT_{n})$. Then, using the above results,
the asymptotic behaviour of the matrix $\rev{\widehat{\mA}}_{\kappa,\mu}^{0}\lbr n\rbr$ is given by
$\rev{\widehat{\mA}}_{\kappa,\mu}^{0}\lbr n\rbr \sim (\mT_{n}^{\#2})^{-1}\cdot\mA_{\kappa,\mu}^{0,\infty}\cdot\mT_{n}^{\#2}$ 
where
\begin{equation*}
\begin{aligned}
\mA_{\kappa,\mu}^{0,\infty} :=
\left\lbr\begin{array}{cc}
0 & \sqrt{\mu/\epsilon}\;\mV_{\kappa}^{0,\infty}\\
\sqrt{\epsilon/\mu}\;\mV_{\kappa}^{0,\infty} & 0
\end{array}\right\rbr =
\left\lbr\begin{array}{cccc}
0 & 0 & 0 & \omega\mu\\
0 & 0 & (\omega\epsilon)^{-1}& 0\\
0 & \omega\epsilon & 0 & 0\\
(\omega\mu)^{-1} & 0 & 0 & 0
\end{array}\right\rbr.
\end{aligned}
\end{equation*}
On the other hand, we also have $\rev{\widehat{\mA}}_{\kappa,\mu}^{1}\lbr n\rbr\sim
(\mT_{n}^{\#2})^{-1}\mA_{\kappa,\mu}^{1,\infty}\mT_{n}^{\#2}$ where
$\mA_{\kappa,\mu}^{1,\infty} := - \mA_{\kappa,\mu}^{0,\infty}$.
Finally, let us define $\mT_{n}^{\#4}\in \mathbb{C}^{8\times 8}$ by $\mT_{n}^{\#4}(\ctru_{1},\ctru_{2}) = (\mT_{n}^{\#2}(\ctru_{1}),\mT_{n}^{\#2}(\ctru_{2}))$
for any $\ctru_{1},\ctru_{2}\in \mathbb{C}^{4}$ i.e. $\mT_{n}^{\#4} = \mrm{diag}(\mT_n^{\#2},\mT_n^{\#2})$. Then we have the asymptotic behaviour 
\begin{equation}
  \begin{aligned}
    & \rev{\widehat{\mrm{MTF}}}_{\loc}\lbr n\rbr\sim (\mT_{n}^{\#4})^{-1}\mrm{MTF}^{\infty}_{\loc}\mT_{n}^{\#4}\quad \text{with}\\
    &\mrm{MTF}_{\loc}^{\infty}:=
    \left\lbr\begin{array}{cc}
    \mA_{\kappa_{0},\mu_{0}}^{0,\infty} & \Id\\[5pt]
    \Id & \mA_{\kappa_{1},\mu_{1}}^{1,\infty}
    \end{array}\right\rbr \;\;\in\C^{8\times 8}
  \end{aligned}
\end{equation}

\begin{remark}
It is important to observe that $\mrm{MTF}_{\loc}^{\infty}$ does
\textit{not} depend on $n$. Since the eigenvalues of
$\rev{\widehat{\mrm{MTF}}}_{\loc}\lbr n\rbr$ coincide with the
eigenvalues of $\mT_{n}^{\#4}\cdot\rev{\widehat{\mrm{MTF}}}_{\loc}\lbr
n\rbr\cdot(\mT_{n}^{\#4})^{-1}$, this shows that the spectrum of
$\rev{\widehat{\mrm{MTF}}}_{\loc}\lbr n\rbr$ converges toward the
spectrum of $\mrm{MTF}_{\loc}^{\infty}$ for $n\to \infty$.
\end{remark}

Now, let us investigate in detail the spectrum of the matrix $\mrm{MTF}_{\loc}^{\infty}$.
First, as an intermediate step, we analyze the spectrum of the matrices:
\begin{align}
\label{eq:def:K}
\rev{\mK^{\infty}_{(\kappa,\mu)}}&:=\mA_{\kappa_{0},\mu_{0}}^{0,\infty}-\mA_{\kappa_{1},\mu_{1}}^{0,\infty},\\
\label{eq:def:S}
\rev{\mS^{\infty}_{(\kappa,\mu)}}&:=\mA_{\kappa_{0},\mu_{0}}^{0,\infty}+\mA_{\kappa_{1},\mu_{1}}^{0,\infty},
\end{align}
where $\rev{\mS^{\infty}_{(\kappa,\mu)}}$ is the Single-Trace Formulation (STF) operator \cite{CAF11}.
A thorough examination shows that they take the form:
\begin{equation}\label{MatrixRumsey}
\mA_{\kappa_{0},\mu_{0}}^{0,\infty}\pm\mA_{\kappa_{1},\mu_{1}}^{0,\infty}=
\left\lbr\begin{array}{cccc}
 0 & 0 & 0 & \alpha^\pm_{\mu}\\
 0 & 0 & \beta^\pm_{\epsilon} & 0\\
 0 & \alpha^\pm_{\epsilon} & 0 & 0\\
 \beta^\pm_{\mu} & 0 & 0 & 0
\end{array}\right\rbr
\quad\text{with}\quad
\left\{\begin{array}{l}
\alpha^\pm_{\epsilon}  = \omega\epsilon_{0} \pm \omega\epsilon_{1}\\
\alpha^\pm_{\mu} = \omega\mu_{0} \pm \omega\mu_{1}\\[5pt]
\beta^\pm_{\epsilon} = (\omega\epsilon_{0})^{-1} \pm(\omega\epsilon_{1})^{-1}\\
\beta^\pm_{\mu} = (\omega\mu_{0})^{-1} \pm (\omega\mu_{1})^{-1}
\end{array}\right.
\end{equation}
Trying to compute directly the eigenvalues of the above matrices leads to the conclusion that any eigenvalue $\lambda^\pm$ satisfies $(\lambda^\pm)^{2}
= \alpha_{\epsilon}^\pm\beta_{\epsilon}^\pm = - (\sqrt{\epsilon_{1}/\epsilon_{0}} \pm \sqrt{\epsilon_{0}/\epsilon_{1}})^{2}$, or 
$(\lambda^\pm)^{2} = \alpha_{\mu}^\pm\beta_{\mu}^\pm=  (\sqrt{\mu_{1}/\mu_{0}} \pm \sqrt{\mu_{0}/\mu_{1}})^{2}$. Thus, the spectr\rev{a} of both matrices is given by
\begin{equation}\label{SpectrumRumsey}
\mathfrak{S}(\rev{\mK^{\infty}_{(\kappa,\mu)}}) = \left\{
\pm\imathj\Lambda_{\mu}, \pm\imathj\Lambda_{\epsilon}\right\}\quad\text{with}\quad
\left\{\begin{array}{l}
\dsp{\Lambda_{\mu}     = \abs{\sqrt{\frac{\mu_{1}}{\mu_{0}}} - \sqrt{\frac{\mu_{0}}{\mu_{1}}}} }\\[15pt]
\dsp{\Lambda_{\epsilon} = \left\vert\sqrt{\frac{\epsilon_{1}}{\epsilon_{0}}}\,-\sqrt{\frac{\epsilon_{0}}{\epsilon_{1}}}\,\right\vert }
\end{array}\right. ,
\end{equation}
\begin{equation}\label{SpectrumSTF}
\mathfrak{S}(\rev{\mS^{\infty}_{(\kappa,\mu)}}) = \left\{
\pm \Upsilon_{\mu}, \pm \Upsilon_{\epsilon}\right\}\quad\text{with}\quad
\left\{\begin{array}{l}
\dsp{\Upsilon_{\mu}     = \abs{\sqrt{\frac{\mu_{1}}{\mu_{0}}} + \sqrt{\frac{\mu_{0}}{\mu_{1}}} }}\\[15pt]
\dsp{\Upsilon_{\epsilon} \,= \abs{\sqrt{\frac{\epsilon_{1}}{\epsilon_{0}}}\,+\sqrt{\frac{\epsilon_{0}}{\epsilon_{1}}}}}
\end{array}\right. .
\end{equation}
Let us now return back to $\mrm{MTF}_{\loc}^{\infty}$. We recall that $(\mA_{\kappa,\mu}^{j,\infty})^{2} =\Id$, and obtain directly the following identity
\begin{equation*}
\left(2\Id-(\mrm{MTF}_{\loc}^{\infty})^{2}\right)^{2} =
\left\lbr\begin{array}{ll}
\rev{(\mK^{\infty}_{(\kappa,\mu)})^{2}} & 0\\
0          & \rev{(\mK^{\infty}_{(\kappa,\mu)})^{2}}
\end{array}\right\rbr.
\end{equation*}
Taking  account of (\ref{SpectrumRumsey}) in addition finally leads to the following expression for the accumulation points of $\mrm{MTF}_{\loc}^\infty$
\begin{equation}\label{eq:eigInfty}
\mathfrak{S}\big(\mrm{MTF}_{\loc}^{\infty}\big) =\left\{\pm\sqrt{2\pm\imathj\Lambda_{\mu} }, \pm\sqrt{2\pm\imathj\Lambda_{\epsilon} }\right\}.
\end{equation}

For numerical experiments and validation, we consider the lossless scattering of Teflon \cite{HAB219} and Ferrite \cite{R87}, both immersed into vacuum. Their relative permeability and permittivity $\epsilon_r,\mu_r$ are described in Fig.~\ref{fig:cases} (left). For each material, we introduce the ``Low'', ``High'' and ``Very High'' frequency regimes with their associate acronyms, corresponding to the excitation of a plane wave with frequency $f$ and wavelength $\lambda:=\frac{2\pi}{\kappa_0}$ as represented in Fig.~\ref{fig:cases} (right).
\begin{figure}[t]
\begin{minipage}{\linewidth}
\begin{minipage}{0.29\linewidth}
\begin{table}[H]
\renewcommand\arraystretch{1.2}
\begin{center}
\footnotesize
\begin{tabular}{|c|c|c|c|} \hline   
  &  $\epsilon_r$ & $\mu_r$ \\ \hline
 Teflon  &  2.1 & 1.0 \\ \hline
 Ferrite  & 2.5 & 1.6\\ \hline
\end{tabular}
\end{center} 
\end{table}  
\end{minipage}
\begin{minipage}{0.7\linewidth}
\begin{table}[H]
\renewcommand\arraystretch{1.2}
\begin{center}
\footnotesize
\begin{tabular}{|c|c|c|c|c|} \hline   
 \multicolumn{2}{|c|}{Case}&  (LF) & (HF) & (VHF) \\ \hline
 \multicolumn{2}{|c|}{$f$} & $50$ MHz & $300$ MHz & $10$ GHz\\ \hline\hline
 \multicolumn{2}{|c|}{$\lambda$ (m)} & $6.0$ & $1.0$ & $0.029$\\ \hline\hline
 \multicolumn{2}{|c|}{$\kappa_0$} & 1.05 & 6.29  & 210\\ \hline 
 \multirow{2}{*}{$\kappa_1$}&Teflon  &  1.52 & 9.11 & 304 \\ \cline{2-4}
 					&Ferrite  & 2.09 & 12.6 & 419\\ \hline
\end{tabular}
\end{center} 
\end{table}  
\end{minipage}
\end{minipage}
\caption{Overview of the material parameters (left) and summary of the frequency $f$, wavelength $\lambda$ and wavenumbers for all cases (right).}
\label{fig:cases}
\end{figure}

We examine numerically the spectrum of the operator $\MTFloc$. An explicit expression 
of the eigenvectors is provided by the vector spherical harmonics $\bfX^{\pp}_{n,m}$ and $\bfX^{\times}_{n,m}$, 
so that $\mathfrak{S}(\MTFloc) = \cup_{n=0}^{+\infty}\mathfrak{S}( \rev{\widehat{\mrm{MTF}}}_{\loc} \lbr n\rbr)$. Each $\mathfrak{S}(\rev{\widehat{\mrm{MTF}}}_{\loc}\lbr n\rbr)$
consists in 8 eigenvalues. On each figure below in Table \ref{tab:Eigenvalue Distribution}, we plot $\cup_{n=0}^{N}\mathfrak{S}(\rev{\widehat{\mrm{MTF}}}_{\loc}\lbr n\rbr)$ (in red) along with the expected accumulation points (in black) for the cases mentioned previously in Fig.~\ref{fig:cases}. We adapt the number of spherical harmonics to the frequency, i.e.~we set $N= 150,200,500$ for the LF, HF and VHF cases, respectively. These plots clearly confirm that: (i) the spectrum has no more than 8 accumulation points that systematically
admit a modulus greater than $\sqrt{2}$; (ii) the accumulation points do not depend on the wavenumber; and (iii) the expected values of the accumulation points coincide with the calculated one. We notice that the eigenvalues spread around the accumulation points and get closer to $0$ with increasing frequency and is likely due to the propagative modes of the local MTF operator (see e.g.~\cite[Section 6]{AD16} for acoustics). The latter induces deterioration of the condition number and of the iteration count for iterative solvers.
\begin{table}[H]
\renewcommand\arraystretch{1.3}
\begin{center}
\footnotesize
\begin{tabular}{
    >{\centering\arraybackslash}m{1cm}
    |>{\centering\arraybackslash}m{6.5cm}
    |>{\centering\arraybackslash}m{6.5cm}
    }
\vspace{0.1cm}
&  Teflon & Ferrite \\ \hline
 (LF) &  \includegraphics[width=1\linewidth]{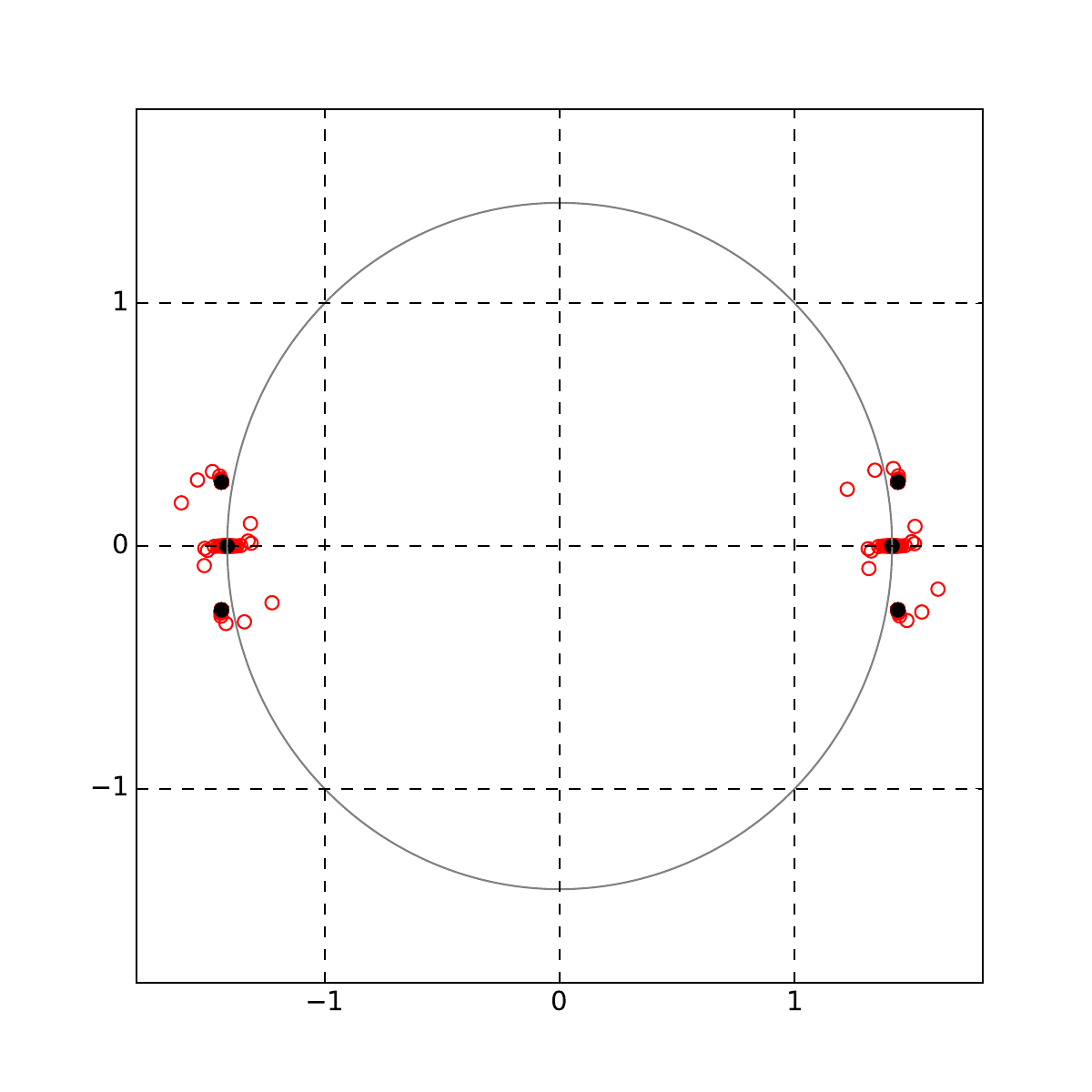} &  \includegraphics[width=1\linewidth]{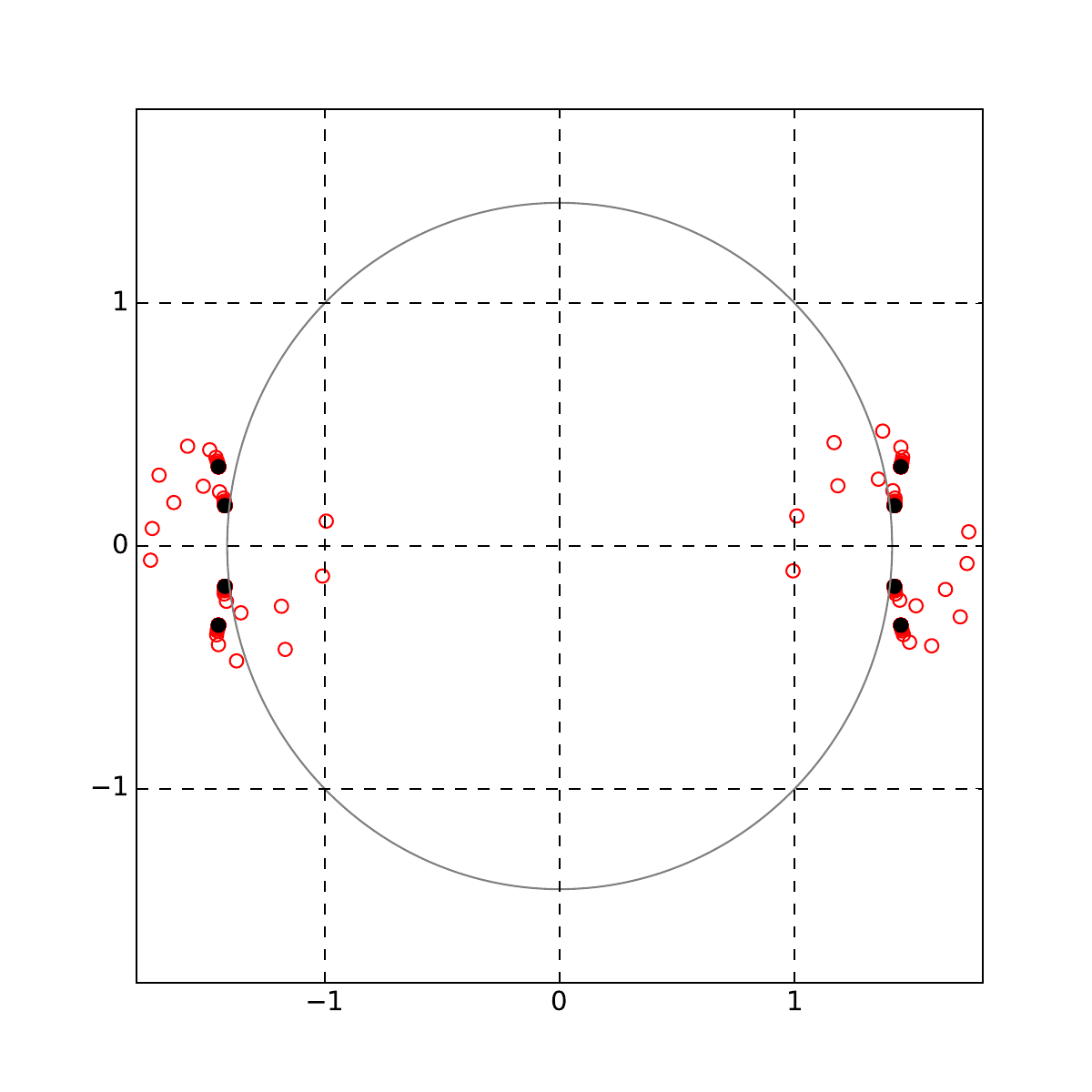}\\ \hline
  (HF) &  \includegraphics[width=1\linewidth]{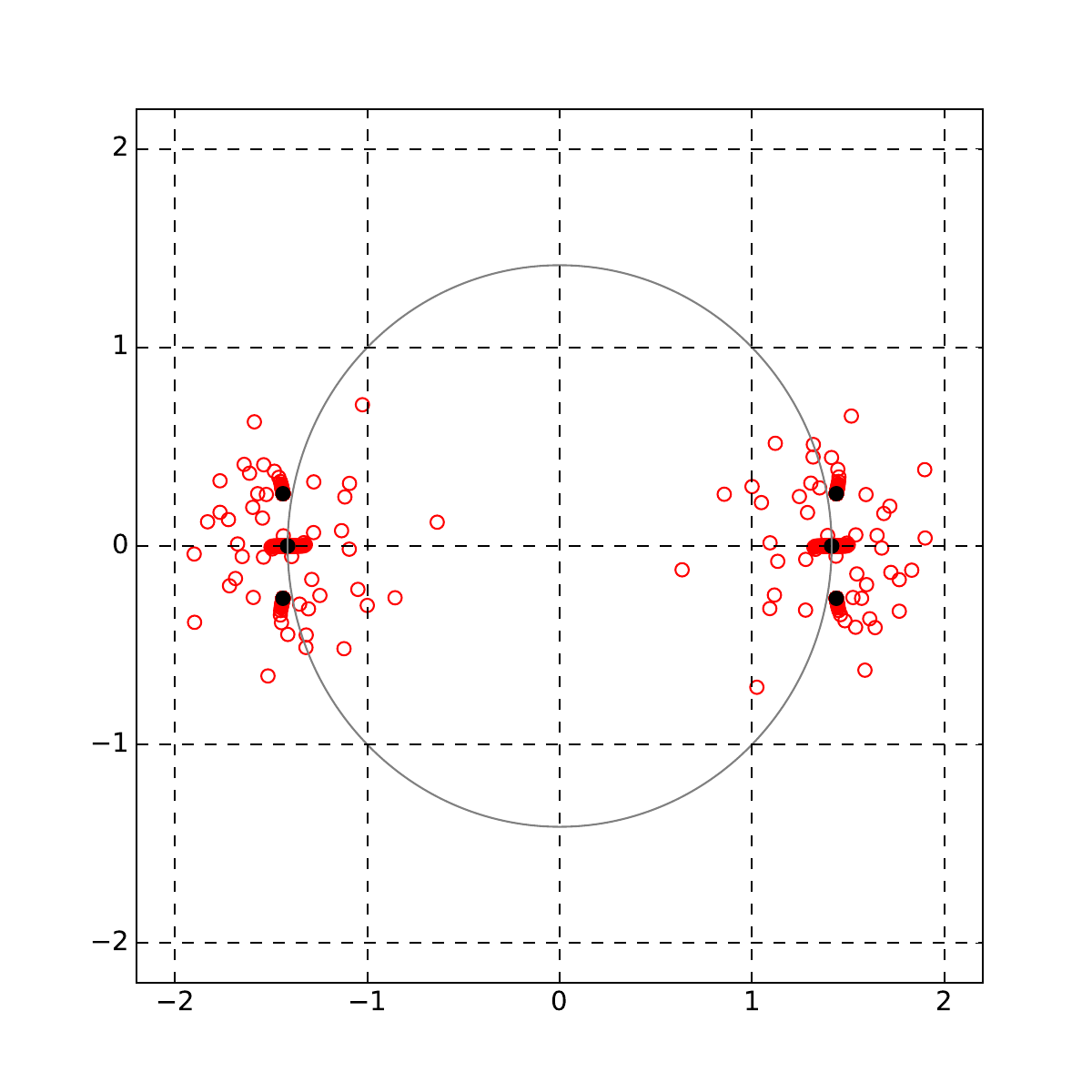} &  \includegraphics[width=1\linewidth]{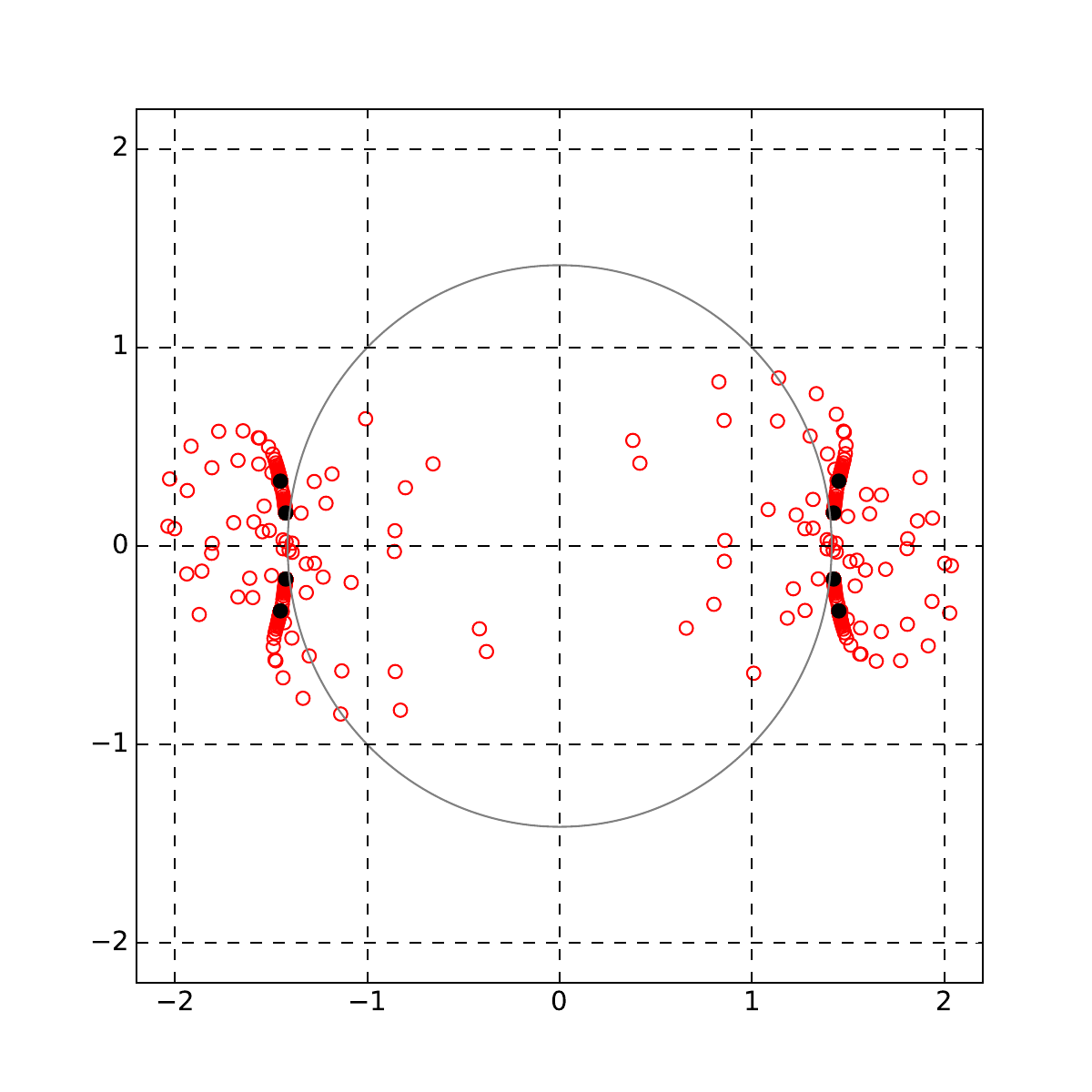}\\ \hline
  (VHF) &  \includegraphics[width=1\linewidth]{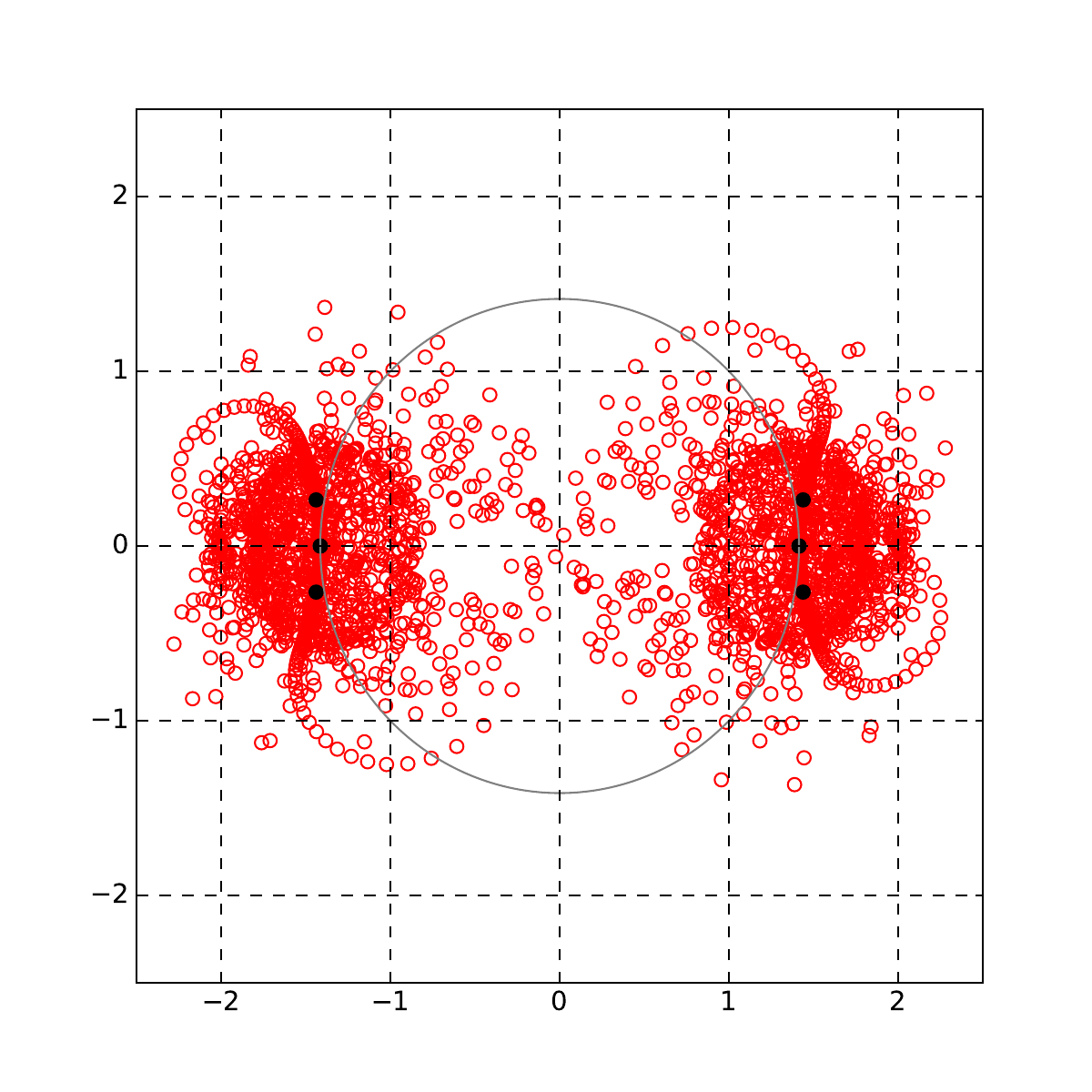} &  \includegraphics[width=1\linewidth]{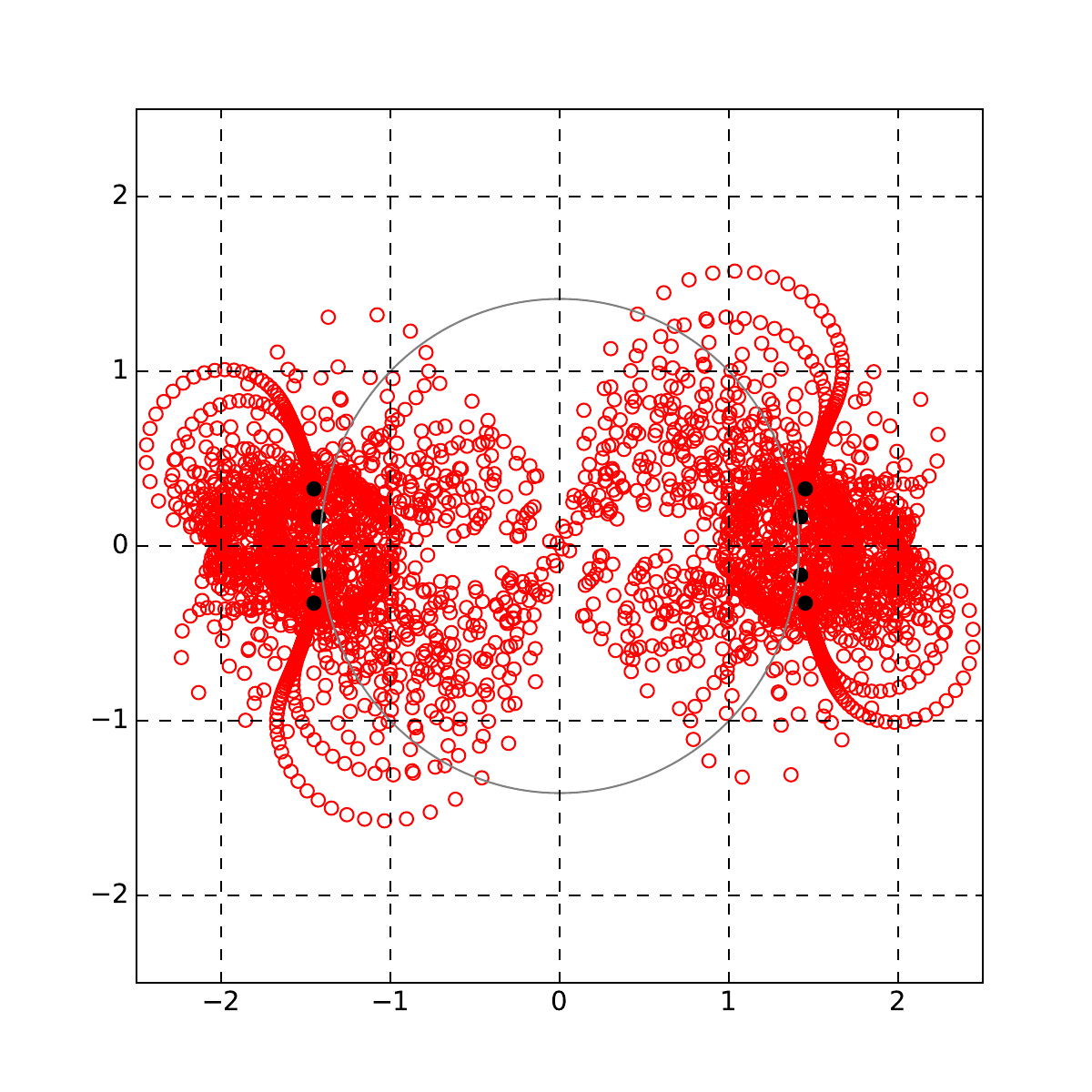}\\ \hline
 \end{tabular}
\caption{Eigenvalue distribution (red) and accumulation points (black) for each case along with a $\sqrt{2}$-radius circle centered at the origin (gray).}
\label{tab:Eigenvalue Distribution}
\end{center} 
\end{table} 

\section{Stability of local MTF for Maxwell equations}\label{StabilityOfLocalMTF}
We now establish a generalized G\aa rding inequality for the local MTF on the unit sphere,
by means of separation of variables. 
First of all, let us derive an expression of the norm on $\bfH^{-\frac{1}{2}}(\div,\Gamma)$ in vector spherical harmonics. Such an expression can be obtained by noting that the dissipative counterpart of the EFIE operator (i.e. associated to a purely imaginary wavenumber) is 
continuous and coercive on $\bfH^{-\frac{1}{2}}(\div,\Gamma)$ so that the corresponding bilinear form 
\begin{equation}
  (\bu,\bv)_{-1/2,\div}:= \int_{\Gamma\times\Gamma}\mathscr{G}_{\imathj}(\bx - \by)\;\big(\div_{\Gamma}\bu(\bx)\div_{\Gamma}
  \overline{\bv}(\by)+\bu(\bx)\cdot\overline{\bv}(\by)\big) d\sigma(\bx,\by)
\end{equation}
yields a scalar product. Here $\mathscr{G}_{\imathj}(\bx) = \exp(-\vert\bx\vert)/(4\pi\vert\bx\vert)$ as $\imathj = \sqrt{-1}$ is
the imaginary unit. The vector fields $\bfX_{n,m}^{\pp}$ and $\bfX_{n,m}^{\times}$ form an orthogonal family with respect to this 
scalar product. As a consequence, to obtain an expression of a norm over  $\bfH^{-\frac{1}{2}}(\div,\Gamma)$, 
one can rely on the decomposition of the dissipative EFIE on vector spherical harmonics. First observe 
that $(\bu,\bv)_{-1/2,\div} = \int_{\Gamma}(\bn_{0}\times\gamma^{0}_{\tan}\cdot\SL_{\kappa}(\bu))\cdot\bv d\sigma $. As a consequence, 
using (\ref{MatrixForm1}) we obtain
\begin{equation}
\begin{aligned}
& (\bu,\bv)_{-1/2,\div}  = \overline{v}^{\top}\cdot\mD_{n}\cdot u\\[5pt]
& \text{where}\quad \mD_{n} = \diag(\mathfrak{J}_{n}'(\imathj)\mathfrak{H}_{n}'(\imathj), \mathfrak{J}_{n}(\imathj)\mathfrak{H}_{n}(\imathj))
\end{aligned}
\end{equation}
\rev{and $u,v \in \C^{2}$ such that  $\bu(\bx) = \bfX_{n,m}(\bx)\cdot u$, and $\bv(\bx) = \bfX_{n,m}(\bx)\cdot v$.}
From this we deduce the asymptotic behaviour $\mD_{n}\sim \widetilde{\mD}_{n}:=\diag( 1+n,1/(1+n))$ for $n\to \infty$, which yields 
the expression of an equivalent norm which is explicit when decomposed in spherical harmonics 
\begin{equation*}
\begin{aligned}
& c_{-}\Vert \bu\Vert_{-1/2,\div}^{2}\leq (\bu,\bu)_{-1/2,\div}\leq c_{+}\Vert \bu\Vert_{-1/2,\div}^{2}\\
& \Vert \bu\Vert_{-1/2,\div}^{2} := \sum_{n=0}^{+\infty}\sum_{\vert m\vert\leq n } \overline{u}_{n,m}^{\top}\cdot\widetilde{\mD}_{n}\cdot u_{n,m}\\
& \text{where}\quad \widetilde{\mD}_{n}:=\diag( 1+n,1/(1+n))
\end{aligned}
\end{equation*}
From this we easily deduce the expression of an explicit norm for $\mbH(\Sigma)$, using the matrix $\mD_{n}^{\#4} := 
\diag(\mD_{n},\mD_{n},\mD_{n},\mD_{n})$. Next we need to introduce intermediate notations for the predominant
behaviour of two key matrices coming into play in the local MTF formulation, namely
\begin{equation}
\begin{aligned}
  \widetilde{\mrm{MTF}}\!\,_{\loc}\lbr n\rbr
  & := (\mT_{n}^{\#4})^{-1}\MTFloc^{\infty}\mT_{n}^{\#4}
  &\hspace{-0.5cm}&\mathop{\sim}_{n\to+\infty}\widehat{\mrm{MTF}}_{\loc}\lbr n\rbr\\
  \widetilde{\mA}_{\kappa_{j},\mu_{j}}^{j}\lbr n\rbr
  &:= (\mT_{n}^{\#2})^{-1}\mA_{\kappa_{j},\mu_{j}}^{j,\infty}\mT_{n}^{\#2}
  &\hspace{-0.5cm}&\mathop{\sim}_{n\to+\infty}\widehat{\mA}_{\kappa_{j},\mu_{j}}^{j}\lbr n\rbr
\end{aligned}
\end{equation}
Since we need to rewrite this formulation variationally, 
we start by inspecting how the duality pairing decomposes on spherical harmonics. First of all, according to 
(\ref{HilbertExpansion}), observe that $\int_{\Gamma}(\bn_{j}\times\bfX^{\pp}_{n,m})\cdot \bfX^{\pp}_{n,m} d\sigma = 
\int_{\Gamma}(\bn_{j}\times\bfX^{\times}_{n,m})\cdot \bfX^{\times}_{n,m} d\sigma = 0$ and $\int_{\Gamma}(\bn_{0}
\times\bfX^{\times}_{n,m})\cdot \bfX^{\pp}_{n,m} d\sigma = 1$. As a consequence, considering the vector fields 
$\ctru(\bx) := \bfX^{\#2}_{n,m}(\bx)\cdot u$ and $\ctrv(\bx) := \bfX^{\#2}_{n,m}(\bx)\cdot v$ where $u,v\in \mathbb{C}^{2}$, 
we have 
\begin{equation}
\begin{aligned}
& \lbr \ctru,\ctrv\rbr_{\Gamma_{0}} =  v^{\top} \mrm{M}\; u,\\
& \text{and}\quad \mrm{M} :=
\left\lbr\begin{array}{rrrr}
0  &  0 &  0  & +1\\
0  &  0 & -1  & 0\\
0  & +1 &  0  & 0\\
-1 &  0 &  0  & 0
\end{array}\right\rbr.
\end{aligned}
\end{equation}
Observe that $\mrm{M}^{\top} = \rev{-}\mrm{M}$. Since $\lbr\ctru,\ctrv\rbr_{\Gamma_{0}} = -\lbr\ctru,\ctrv\rbr_{\Gamma_{1}}$, 
we obtain a global matrix expression for the pairing on the multi-trace space: for  
$\ctru(\bx) := \bfX^{\#4}_{n,m}(\bx)\cdot u$ and $\ctrv(\bx) := \bfX^{\#4}_{n,m}(\bx)\cdot v$ where $u,v\in \mathbb{C}^{4}$, 
we have 
\begin{equation}
\llbr \ctru,\ctrv\rrbr_{\Sigma} =  v^{\top} \mathbb{M}\; u
\quad \text{and}\quad \mathbb{M} :=
\left\lbr\begin{array}{cc}
+\mrm{M} & 0 \\ 0 & -\mrm{M}
\end{array}\right\rbr.
\end{equation}
To examine coercivity of local MTF on the sphere, we need to study  the coercivity of the matrix 
$\mathbb{M}\cdot\widehat{\mrm{MTF}}_{\loc}\lbr n\rbr$ as $n\to\infty$. If we look at the asymptotic behaviour of this 
matrix, taking account of the results of Section \ref{Asymptotics}, we obtain the expression 
\begin{equation}\label{ExprDiagBlocks}
(-1)^{j}\,\mrm{M}\cdot\widetilde{\mA}_{\kappa_{j},\mu_{j}}^{j}\lbr n\rbr\; =  
\left\lbr\begin{array}{cccc}
\dsp{\frac{n}{\omega\mu_{j}}}  &  0                            &  0                             & 0\\
0                       &  -\dsp{\frac{\omega\epsilon_{j}}{n}} &  0                             & 0\\
0                       &  0                                  &  \dsp{\frac{n}{\omega\epsilon_{j}}}  & 0\\
0                       &  0                                  &  0                             & -\dsp{\frac{\omega\mu_{j}}{n}}
\end{array}\right\rbr.
\end{equation}
Let us introduce a diagonal matrix $\theta\in\mathbb{R}^{2\times 2}$ defined by $\theta = \mrm{diag(+1,-1)}$, 
and denote $\Theta := \rev{\mrm{diag}(\theta,\theta)}\in\mathbb{R}^{4\times 4}$. From (\ref{ExprDiagBlocks}), it clearly follows 
that  $(-1)^{j}\,\mrm{M}\cdot\widetilde{\mA}_{\kappa_{j},\mu_{j}}^{j}\lbr n\rbr\cdot \Theta$ is a real valued diagonal positive definite matrix. 
On the other hand $(\mathbb{M}\cdot\Theta)^{\top} = - \mathbb{M}\cdot\Theta$. As a consequence we finally conclude
that there exists $c>0$ independent of $n$ such that 
\begin{equation}
\begin{aligned}
& \Re e\{\mrm{U}^{\top}\cdot\widetilde{\mrm{MTF}}_{\loc}\lbr n\rbr\cdot\Theta \cdot \overline{\mrm{U}}\}\geq c\;
\mrm{U}^{\top}\cdot\mD^{\#4}_{n}\cdot \overline{\mrm{U}}\\
&  \forall \mrm{U}\in\mathbb{C}^{8},\;\forall n\geq 0.
\end{aligned}
\end{equation}
Since the constant $c>0$ is independent of $n$, summing this inequality over $n$, and taking account that
$\widetilde{\mrm{MTF}}\!\,_{\loc}^{\infty}\lbr n\rbr$ is the asymptotic behaviour of $\widehat{\mrm{MTF}}\!\,_{\loc}^{\infty}\lbr n\rbr$, we finally obtain the following coercivity statement.   

\begin{theorem}\label{GardingIneq}
  There exists a compact operator $\mathcal{K}:\mbH(\Sigma)\to\mbH(\Sigma)$
  and a constant $C>0$ such that for all $\ctru\in \mbH(\Sigma)$ we have
  \begin{equation}
    \Re e\left\{\llbr (\MTFloc+\mathcal{K})\bctru,\Theta(\overline{\bctru})\rrbr\right\}\geq C\Vert \bctru\Vert_{\mbH(\Sigma)}^{2}.
  \end{equation}
\end{theorem}

\section{Preconditioning the local MTF for Maxwell equations}\label{sec:Preconditioning}
In this section, we introduce a closed formula for the inverse of the multi-trace operator and propose robust preconditioners for the formulation. First, let us rewrite $\mrm{MTF}_{\loc}$ as:
\begin{align}\quad \mrm{MTF}_{\loc}= 
\left\lbr\begin{array}{cc}
\mA_{\kappa_{0},\mu_{0}}^{0} & \Id\\[5pt]
\Id & -\mA_{\kappa_{1},\mu_{1}}^{0}
\end{array}\right\rbr,
\end{align}
and introduce the block diagonal STF operator of (\ref{eq:def:S})
\begin{equation}
\IS_{(\kappa,\mu)} := \mrm{diag}(\rev{\mS_{(\kappa,\mu)}},\rev{\mS_{(\kappa,\mu)}}) \text{ with } \rev{\mS_{(\kappa,\mu)}} :=(\mA_{\kappa_{0},\mu_{0}}^{0}+\mA_{\kappa_{1},\mu_{1}}^{0}),
\end{equation}
with $\rev{\mS_{(\kappa,\mu)}}$ known to be invertible and $(\rev{\mS_{(\kappa,\mu)}})^2$ a second-kind Fredholm operator for smooth surfaces.
Finally, we also introduce $\IK_{(\kappa,\mu)}:=\mrm{diag}(\rev{\mK_{(\kappa,\mu)}},\rev{\mK_{(\kappa,\mu)}})$ with $\rev{\mK_{(\kappa,\mu)}} := (\mA_{\kappa_{0},\mu_{0}}^{0}-\mA_{\kappa_{1},\mu_{1}}^{0})$ \rev{in (\ref{eq:def:K}) being} a compact operator on smooth surfaces.

\rev{We recall the inverse formula for $2\times 2$ matrices, combining (i) and (ii) in \cite[Theorem 2.1]{lu2002inverses}, and being also valid for bounded linear operators.}
\rev{\begin{lemma}[{\cite[Theorem 2.1]{lu2002inverses}}]\label{LemmaInverse}Set 
$$
M= \left\lbr\begin{array}{cc}
A & B\\[5pt]
C &D
\end{array}\right\rbr
$$
with $A$,$D$,$D- CA^{-1}B$ and $A - BD^{-1}C$ being invertible. Then $M$ is invertible and
$$
M^{-1}= \left\lbr\begin{array}{cc}
(A - BD^{-1}C)^{-1} & - A^{-1} B (D- CA^{-1}B)^{-1}\\[5pt]
-D^{-1}C(A-BD^{-1}C)^{-1} & (D-CA^{-1}B)^{-1}
\end{array}\right\rbr.
$$
\end{lemma}}

\begin{theorem}
\label{ClosedInverse} 
The exact inverse of the multi-trace operator is given by 
\begin{align}\label{MTFLocInv}
\mrm{MTF}_{\loc}^{-1}= \left\lbr\begin{array}{cc}
\rev{\mS_{(\kappa,\mu)}^{-1}} & \mA_{\kappa_{0},\mu_{0}}^{0}\rev{\mS_{(\kappa,\mu)}^{-1}}\\[5pt]
\mA_{\kappa_{1},\mu_{1}}^{0}\rev{\mS_{(\kappa,\mu)}^{-1}} & -\rev{\mS_{(\kappa,\mu)}^{-1}}
\end{array}\right\rbr.
\end{align}
\end{theorem}
\begin{proof}
\rev{This result is straightforward by application of Lemma \ref{LemmaInverse} to
\begin{align}\quad \mrm{MTF}_{\loc}= 
\left\lbr\begin{array}{cc}
\mA_{\kappa_{0},\mu_{0}}^{0} & \Id\\[5pt]
\Id & -\mA_{\kappa_{1},\mu_{1}}^{0}
\end{array}\right\rbr,
\end{align}
and using identities $(\mA^0_{\kappa_0,\mu_0})^2=(\mA^0_{\kappa_1,\mu_1})^2=\Id$, leading to
\begin{align*}
\mrm{MTF}_{\loc}^{-1}
&= \left\lbr\begin{array}{cc}
(\mA^0_{\kappa_0,\mu_0}+\mA^0_{\kappa_1,\mu_1})^{-1} & \mA^0_{\kappa_0,\mu_0}(\mA^0_{\kappa_1,\mu_1}+\mA^0_{\kappa_0,\mu_0})^{-1}\\[5pt]
\mA^0_{\kappa_1,\mu_1}(\mA^0_{\kappa_0,\mu_0}+\mA^0_{\kappa_1,\mu_1})^{-1} & -(\mA^0_{\kappa_1,\mu_1}+\mA^0_{\kappa_0,\mu_0})^{-1}
\end{array}\right\rbr\\
&= \left\lbr\begin{array}{cc}
\mS_{(\kappa,\mu)}^{-1} & \mA_{\kappa_{0},\mu_{0}}^{0}\mS_{(\kappa,\mu)}^{-1}\\[5pt]
\mA_{\kappa_{1},\mu_{1}}^{0}\mS_{(\kappa,\mu)}^{-1} & -\mS_{(\kappa,\mu)}^{-1}
\end{array}\right\rbr.
\end{align*}}
\end{proof}
From \rev{Theorem \ref{ClosedInverse}} we learn that the multi-trace operator can be closely related to the
inverse of the single-trace operator. Now, remembering that \rev{$\mS^2_{(\kappa,\mu)}$} is a compact perturbation of the
identity for smooth surfaces, it could be appropriate to replace \rev{$\mS^{-1}_{(\kappa,\mu)}$} by \rev{$\mS_{(\kappa,\mu)}$}. Using Theorem \ref{ClosedInverse}, we state the following important result:
\begin{proposition}
Introduce the following operator:
\begin{align}\label{Prec}
\rev{\IB_{(\kappa,\mu)}}:= \left\lbr\begin{array}{cc}
\rev{\mS_{(\kappa,\mu)}} & \mA_{\kappa_{0},\mu_{0}}^{0}\rev{\mS_{(\kappa,\mu)}}\\[5pt]
\mA_{\kappa_{1},\mu_{1}}^{0}\rev{\mS_{(\kappa,\mu)}} & -\rev{\mS_{(\kappa,\mu)}}
\end{array}\right\rbr.
\end{align}
Then $\rev{\IB_{(\kappa,\mu)}} \cdot\mrm{MTF}_{\loc}= \rev{\IS}_{(\kappa,\mu)}^2$. Also, if $\kappa_0=\kappa_1$ and $\mu_0 = \mu_1$,
then $\rev{\IB_{(\kappa,\mu)}} \cdot \mrm{MTF}_{\loc}   = 2 \Id$.
\end{proposition}

In parallel, we introduce the usual squared operator preconditioner and detail its properties:
\begin{theorem}
The square of $\mrm{MTF}_{\loc}$ is given by:
\begin{align}
\mrm{MTF}_{\loc}^{2}= \left\lbr\begin{array}{cc}
2 \Id &  \rev{\mK_{(\kappa,\mu)}}\\[5pt]
\rev{\mK_{(\kappa,\mu)}} & 2 \Id
\end{array}\right\rbr= 2\Id + \Pi \rev{\IK}_{(\kappa,\mu)}
\end{align}
\end{theorem}
\noindent 
In addition, in \cite[\S 5.1]{MR3614903} it was already pointed that,
if $\kappa_0=\kappa_1$ and $\mu_0 = \mu_1$, then $ \mrm{MTF}^2_{\loc} = 2\Id $.
In \cite{CLAEYS2020} such properties were used to investigate the close
relationship between local MTF and optimized Schwarz methods. 
\begin{remark}
\label{rmk:block}
Sometimes, the MTF is preconditioned by $\rev{\IA}_{(\kappa,\mu)}$, which provides
\begin{align}
\rev{\IA}_{(\kappa,\mu)} \cdot \mrm{MTF}_{\loc} = \left\lbr\begin{array}{cc}
 \Id &  \mA_{\kappa_{0},\mu_{0}}^{0}\\[5pt]
\mA_{\kappa_{1},\mu_{1}}^{0} &  \Id
\end{array}\right\rbr= \Id + \rev{\IA}_{(\kappa,\mu)} \Pi ,
\end{align}
Similarly, one could use $\Pi$ as a preconditioner, which is a cost-free alternative due to the sparse nature of this operator. Those solutions allow to obtain accumulation points with positive real part but are not a compact perturbation of identity. We did not incorporate these preconditioner in our analysis due to GMRes iteration counts close to $\mrm{MTF}_{\loc}$.
\end{remark}

\subsection{Preconditioning: Clustering properties}

The novel preconditioner \rev{$\IB_{(\kappa,\mu)}$} proposed appears to be a second-order approximation of the inverse operator while $\mrm{MTF}_{\loc}^{2}$ could be considered as a first-order approximation of the inverse operator:
\begin{align*}
\rev{\IB_{(\kappa,\mu)}} \cdot \mrm{MTF}_{\loc}&= 2 \Id + \rev{\IK_1^2}\textup{, and}\\
\mrm{MTF}_{\loc}^2 & = 2 \Id + \rev{\IK_2},
\end{align*}
with $\rev{\IK_1,\IK_2}:\bfH^{-\frac{1}{2}}(\div,\Gamma)^{2}\to \bfH^{-\frac{1}{2}}(\div,\Gamma)^{2}$ compact operators. We can expect this second-order property---relatively to a first-order property---to imply: (i) faster convergence towards zero of the singular values that lay close to the cluster and (ii) increasing spreading of the outlying singular values, with direct consequences on iterative solvers. \rev{We refer to \cite[Section 5]{EIJH21} for results concerning the convergence of iterative solvers applied to (operator) preconditioned schemes.}

Notice that $\rev{\IB_{(\kappa,\mu)}}$ does not involve new operators and is straightforwardly computable from the knowledge of $\mrm{MTF}_{\loc}$. Still, it involves another operator product, which would originate a preconditioner that consists in two matrix-vector products in case of discretization with adapted function spaces. Taking  account of (\ref{SpectrumRumsey}) finally leads to the following expression for the accumulation points of all aforementioned operators
\begin{equation}\label{eq:accPoints}
\begin{array}{rcl}
\mathfrak{S}\big(\mrm{MTF}_{\loc}^{\infty}\big) &= &\{\pm\sqrt{2\pm\imathj\Lambda_{\mu} }, \pm\sqrt{2\pm\imathj\Lambda_{\epsilon} }\},\\
\mathfrak{S}\big((\mrm{MTF}_{\loc}^{\infty})^2\big) &=&\{2\pm\imathj\Lambda_{\mu} , 2\pm\imathj\Lambda_{\epsilon}\},\\
\mathfrak{S}\big(\rev{\IB_{(\kappa,\mu)}} \cdot \mrm{MTF}_{\loc}^{\infty}\big) &=&\{
\Upsilon_{\mu}^2,\Upsilon_{\epsilon}^2\}.
\end{array}
\end{equation}
Accumulation points for the last formulation are surprisingly simple as they rewrite as
\begin{align*}
\left\{\Upsilon_{\mu}^2,\Upsilon_{\epsilon}^2\right\} & = \left\{ 2 + \mu_r + 1/\mu_r, 2 + \epsilon_r + 1/\epsilon_r\right\}.
\end{align*}
Besides, we define $\underline{\Upsilon} :=  \min  (\Upsilon_{\mu}, \Upsilon_{\epsilon})>1$. As stated before, the local MTF operator has no more than eight accumulation points, the latter being reduced to $4$ accumulation points when using squared operator preconditioning, under the requirement of performing a two matrix-vector products at each iteration of iterative solvers. Their barycenter is located at $2.0$ independently of the medium parameters, allowing for further clustering properties (see \cite{CIK96} for the analogy between one big cluster and several small ones). Finally, the novel approximate inverse has not more than two accumulation points, whose center is bounded away from zero, at the price of performing two additional matrix-vector products at each iteration of linear solvers.
Notice that the two accumulation points of $\rev{\IB_{(\kappa,\mu)} }\cdot \mrm{MTF}_{\loc}$  and their midpoints are parameter
dependent. Still, these can be rescaled by a factor $\underline{\Upsilon}^{2}$ if needed to bring their
values closer to one.

\begin{figure}[H]
  \centering
  \includegraphics[width=1\linewidth,height=6cm]{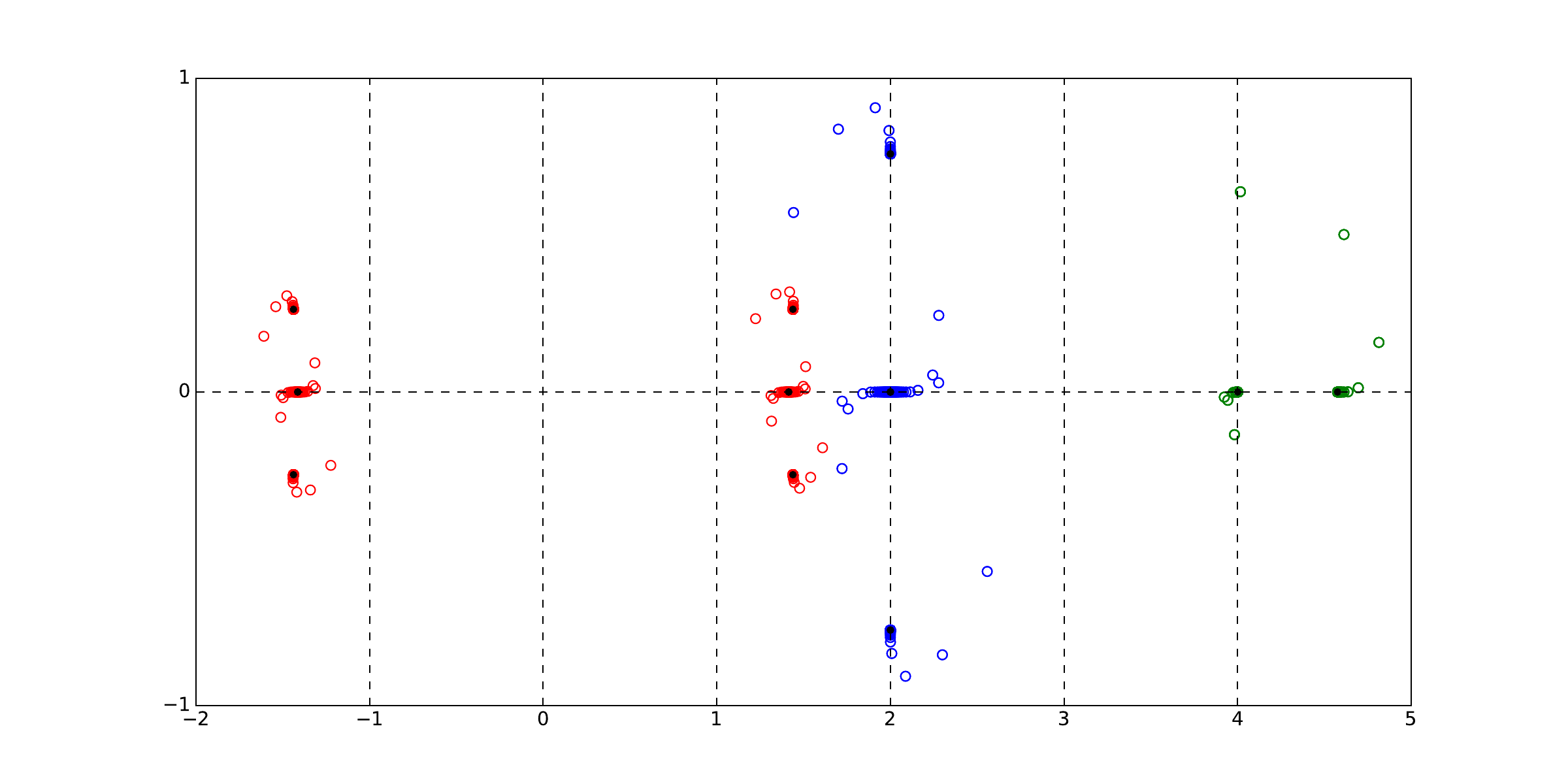}
\caption{Eigenvalue distribution for the Teflon LF case of $\mrm{MTF}_{\loc}$ (red), $\mrm{MTF}_{\loc}^2$ (blue) and $\rev{\IB_{(\kappa,\mu)}} \cdot \mrm{MTF}_{\loc}$ (green).}
\label{fig:eigenvalues_prec}
\end{figure}

In Fig.~\ref{fig:eigenvalues_prec}, we plot the eigenvalue distribution for the preconditioned operators for the Teflon and the LF case and remark that the accumulation points coincide with their expected values.

Next, we decide to normalize the matrices for comparison purposes, namely, we compare $\sqrt{2}^{-1} \mrm{MTF}_{\loc}$, $\mrm{MTF}_{\loc}^2 /2$ and $\underline{\Upsilon}^{-2}\rev{\IB_{(\kappa,\mu)}} \cdot \mrm{MTF}_{\loc}$. In Table.~\ref{tab:EigenvalueDistributionPrec}, we represent the eigenvalue distribution of all proposed operators for the three frequency ranges. These are significant, and show that from left to right: (i) the number of accumulation points diminishes and we observe stronger clustering close to the accumulation points while (ii) the outlying eigenvalues are more spread for increasing frequencies as expected. 

\begin{table}[H]
\renewcommand\arraystretch{1.4}
\begin{center}
\footnotesize
\begin{tabular}{
    >{\centering\arraybackslash}m{1cm}|
    >{\centering\arraybackslash}m{4.5cm}|
    >{\centering\arraybackslash}m{4.5cm}|
    >{\centering\arraybackslash}m{4.5cm}
    }
\vspace{0.1cm}
  & $\mrm{MTF}_{\loc}$ & $(\mrm{MTF}_{\loc})^2$ & $\rev{\IB_{(\kappa,\mu)}} \cdot \mrm{MTF}_{\loc}$\\ \hline
 (LF) &  \includegraphics[width=1\linewidth]{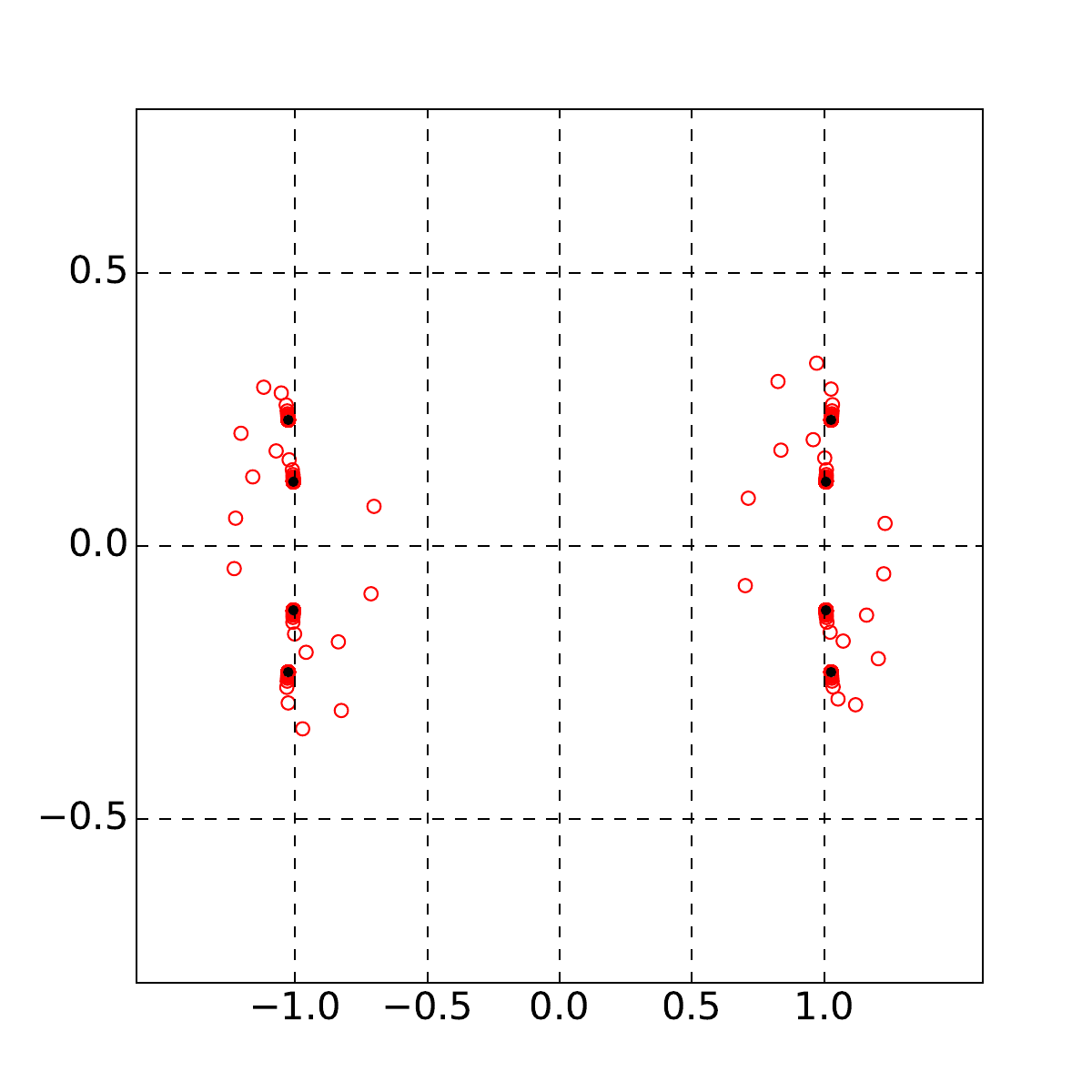}& \includegraphics[width=1\linewidth]{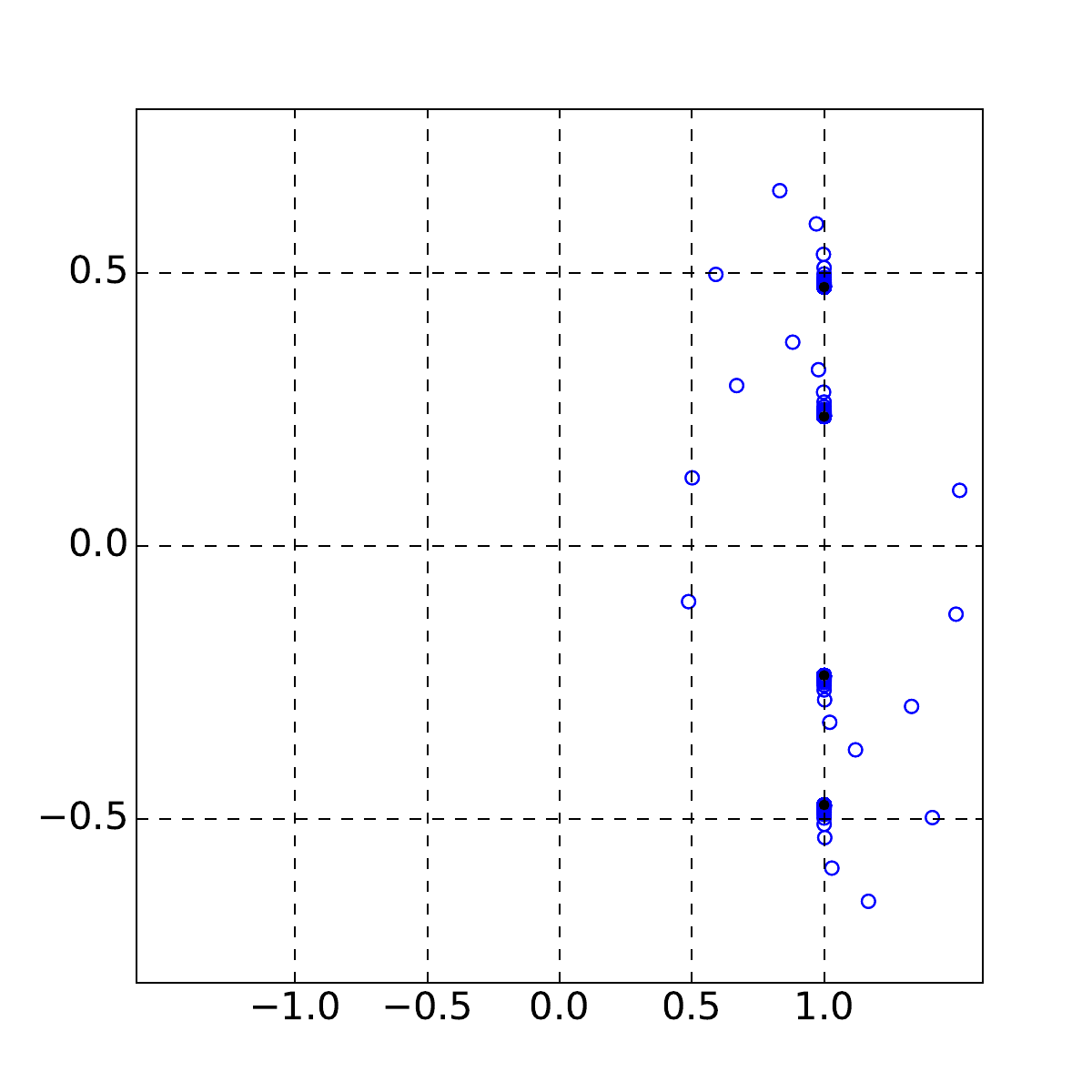}& \includegraphics[width=1\linewidth]{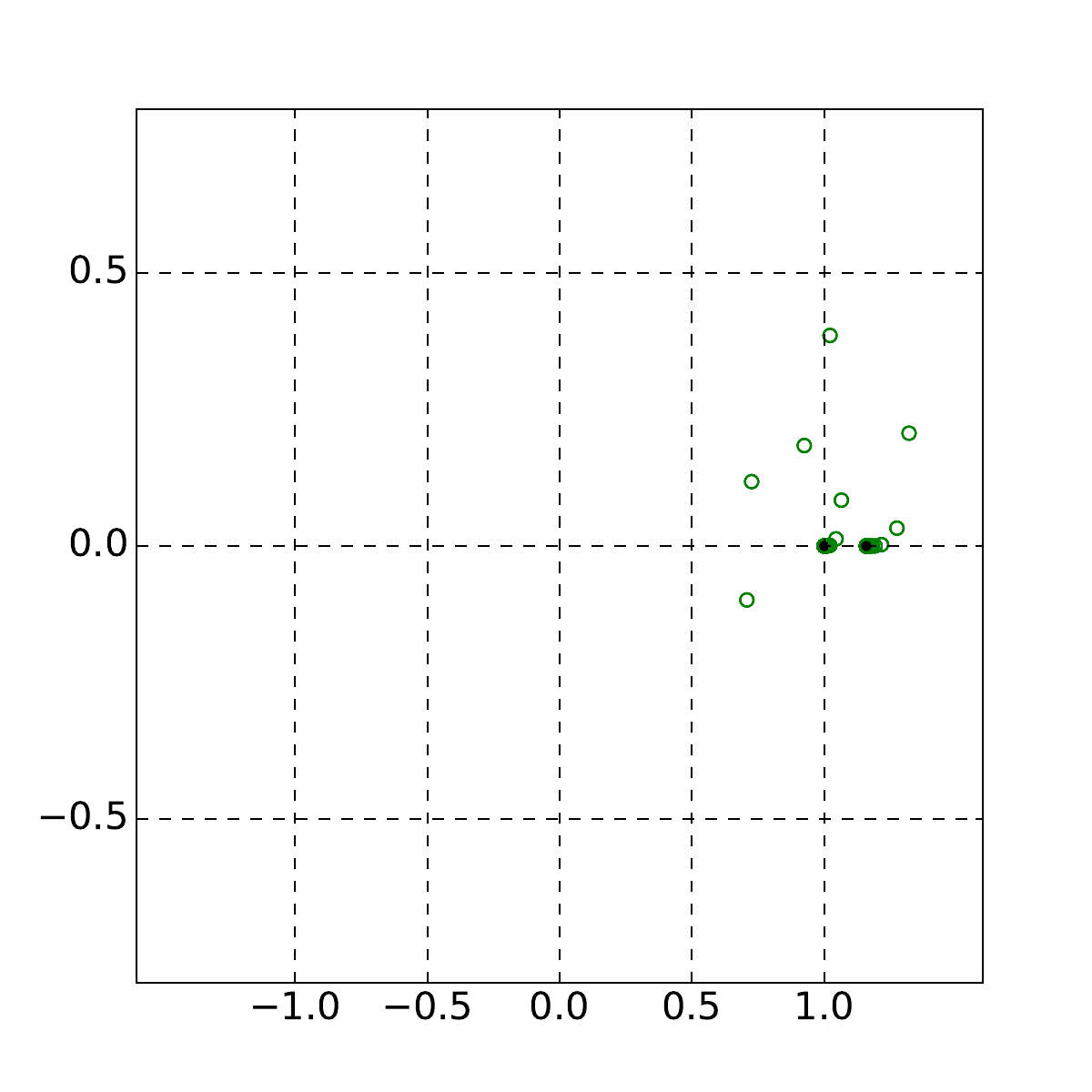}\\ \hline
 (HF) &  \includegraphics[width=1\linewidth]{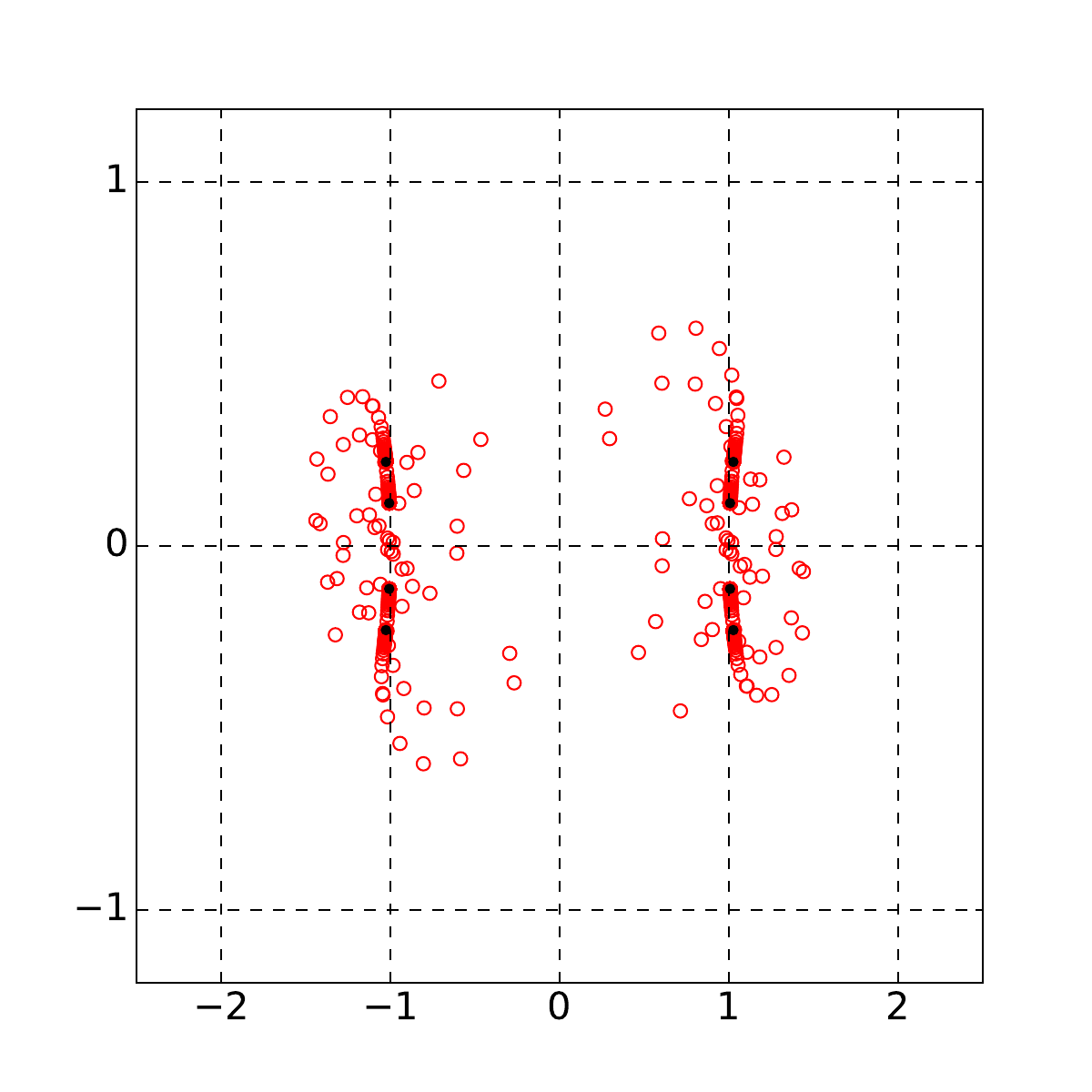}& \includegraphics[width=1\linewidth]{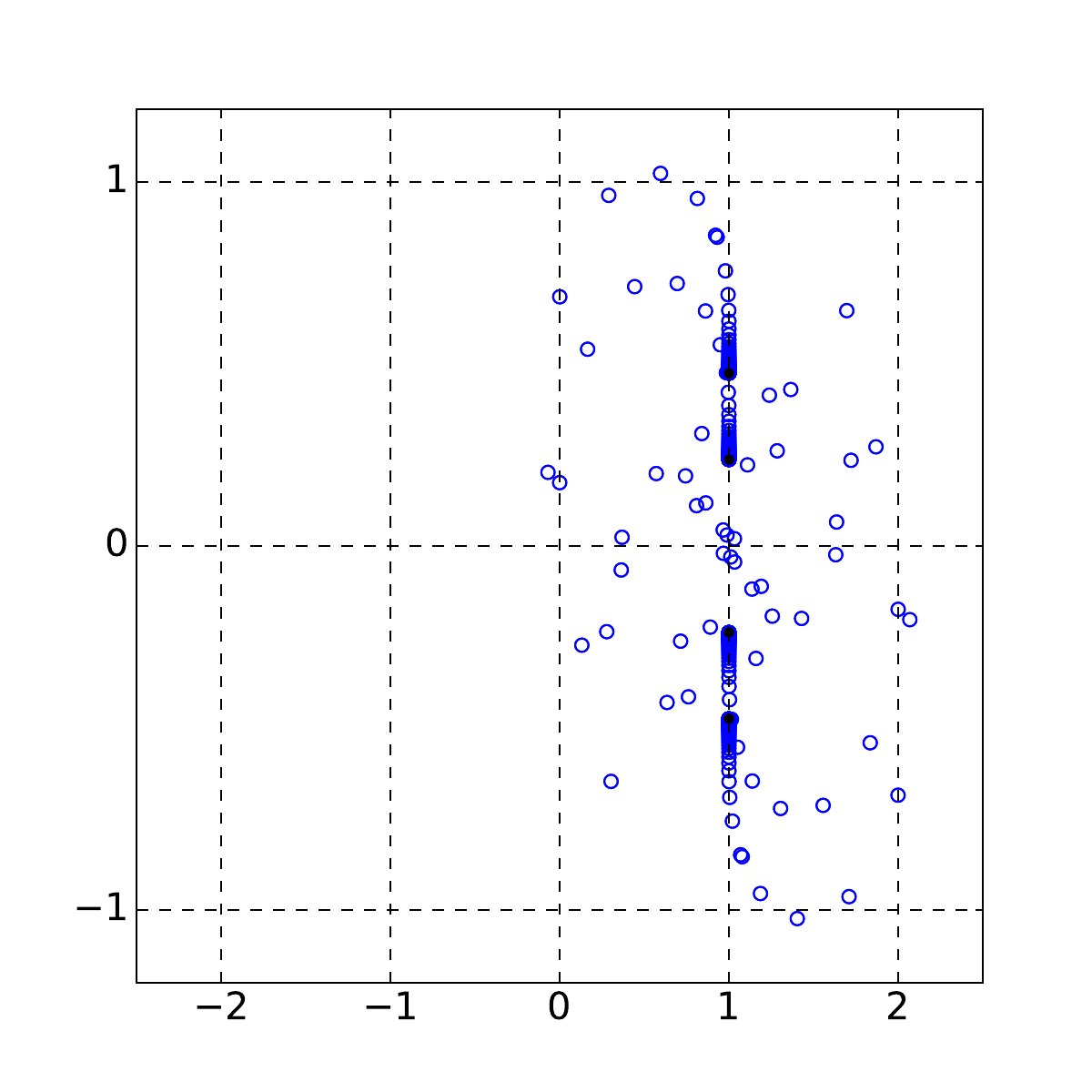}& \includegraphics[width=1\linewidth]{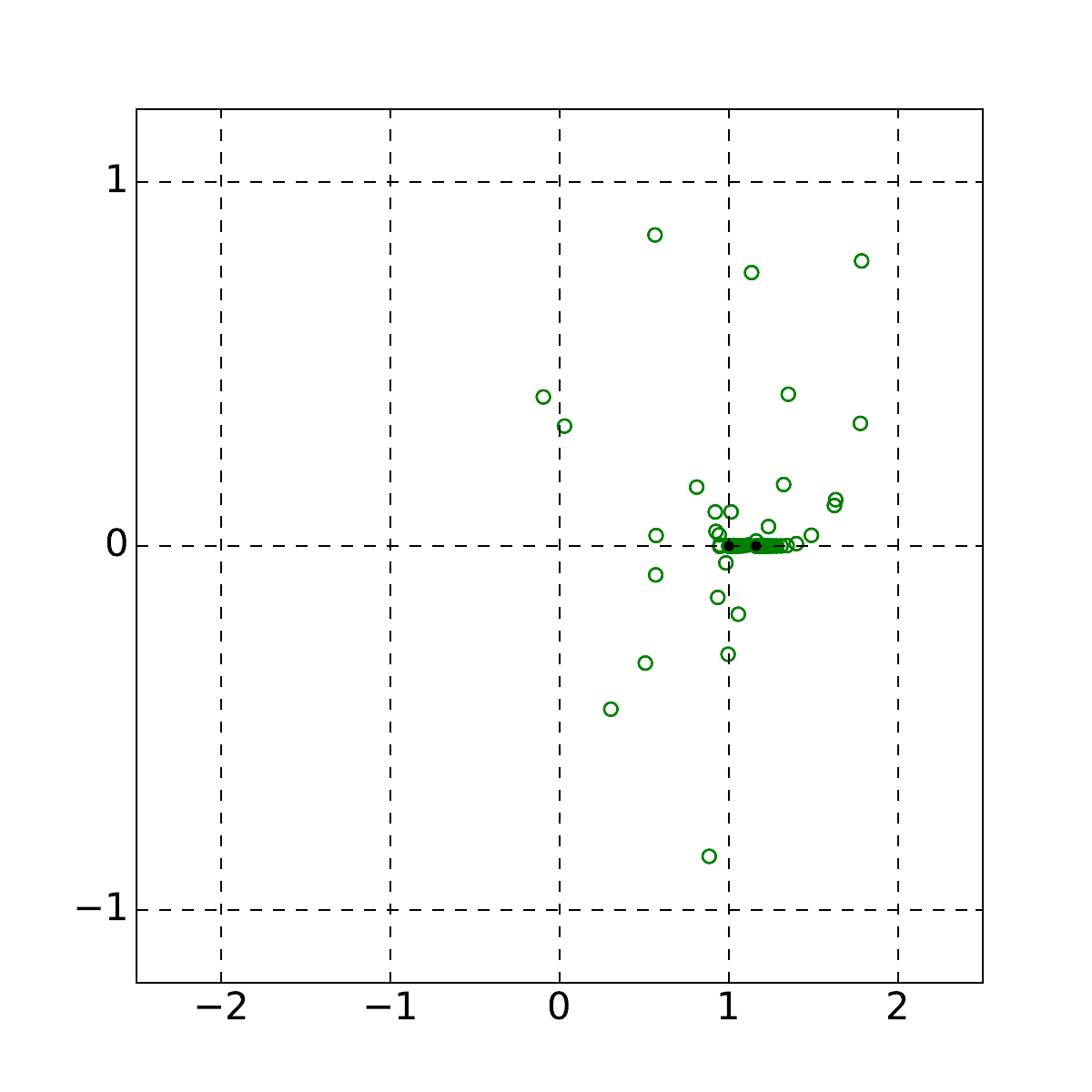}\\ \hline
 (VHF) &  \includegraphics[width=1\linewidth]{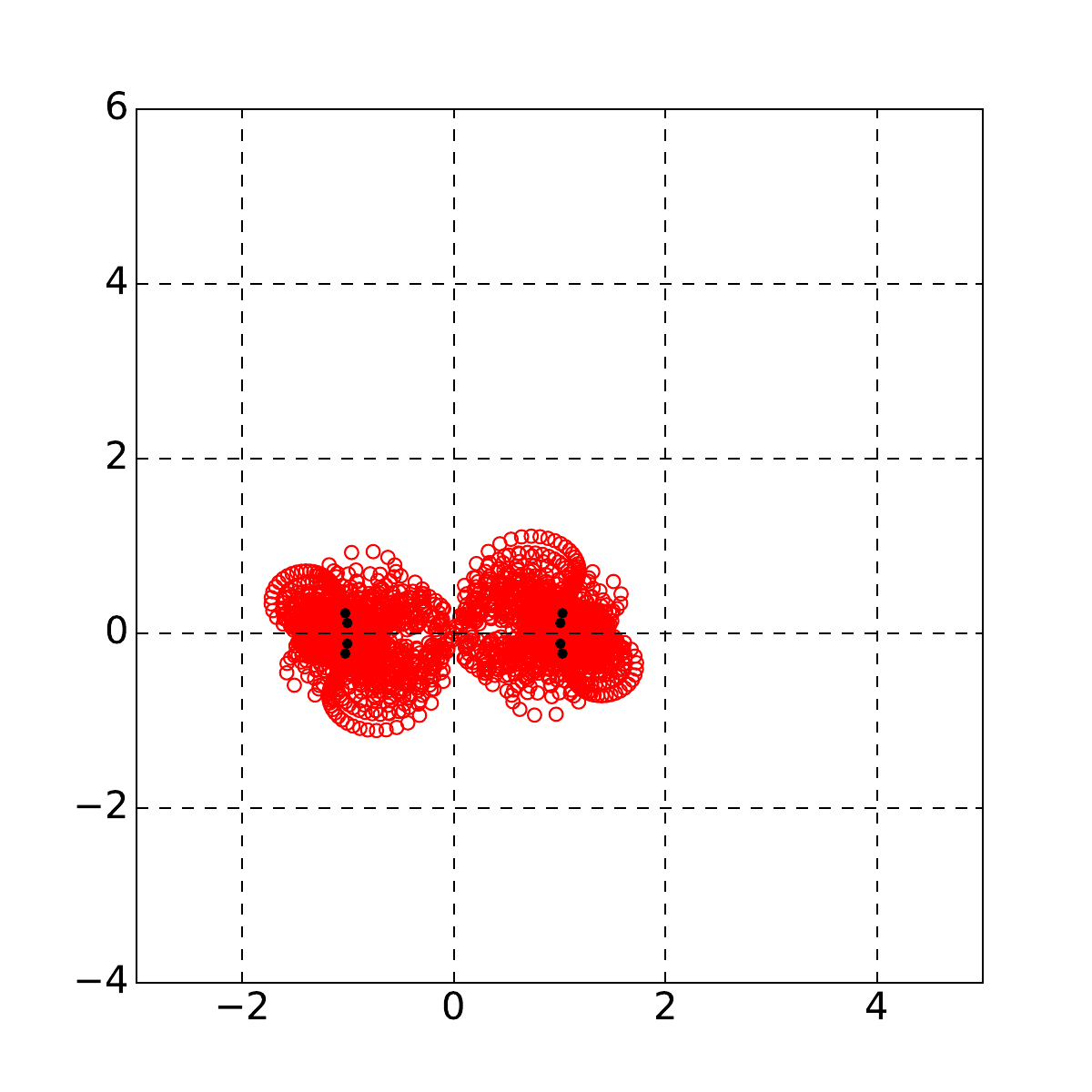}& \includegraphics[width=1\linewidth]{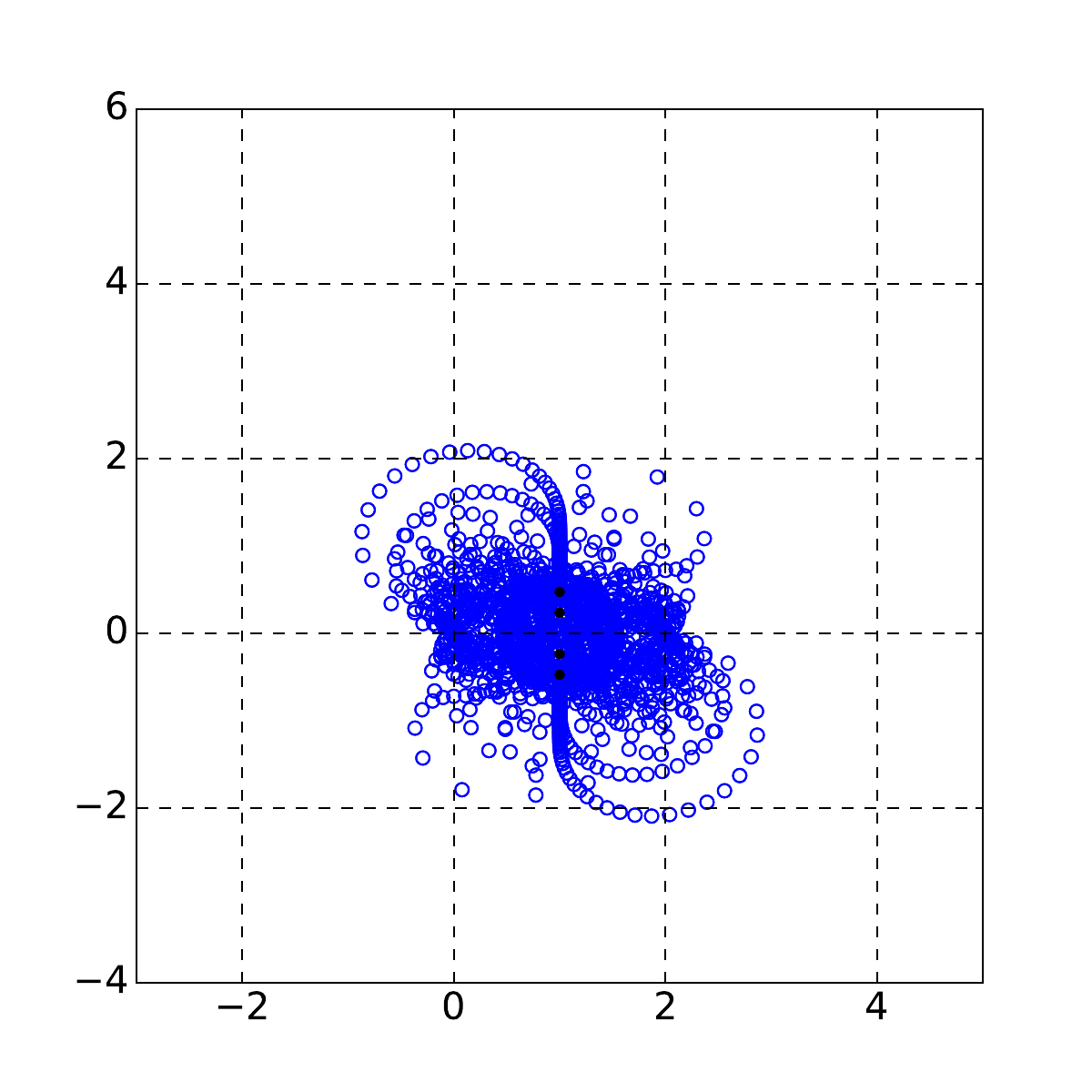}& \includegraphics[width=1\linewidth]{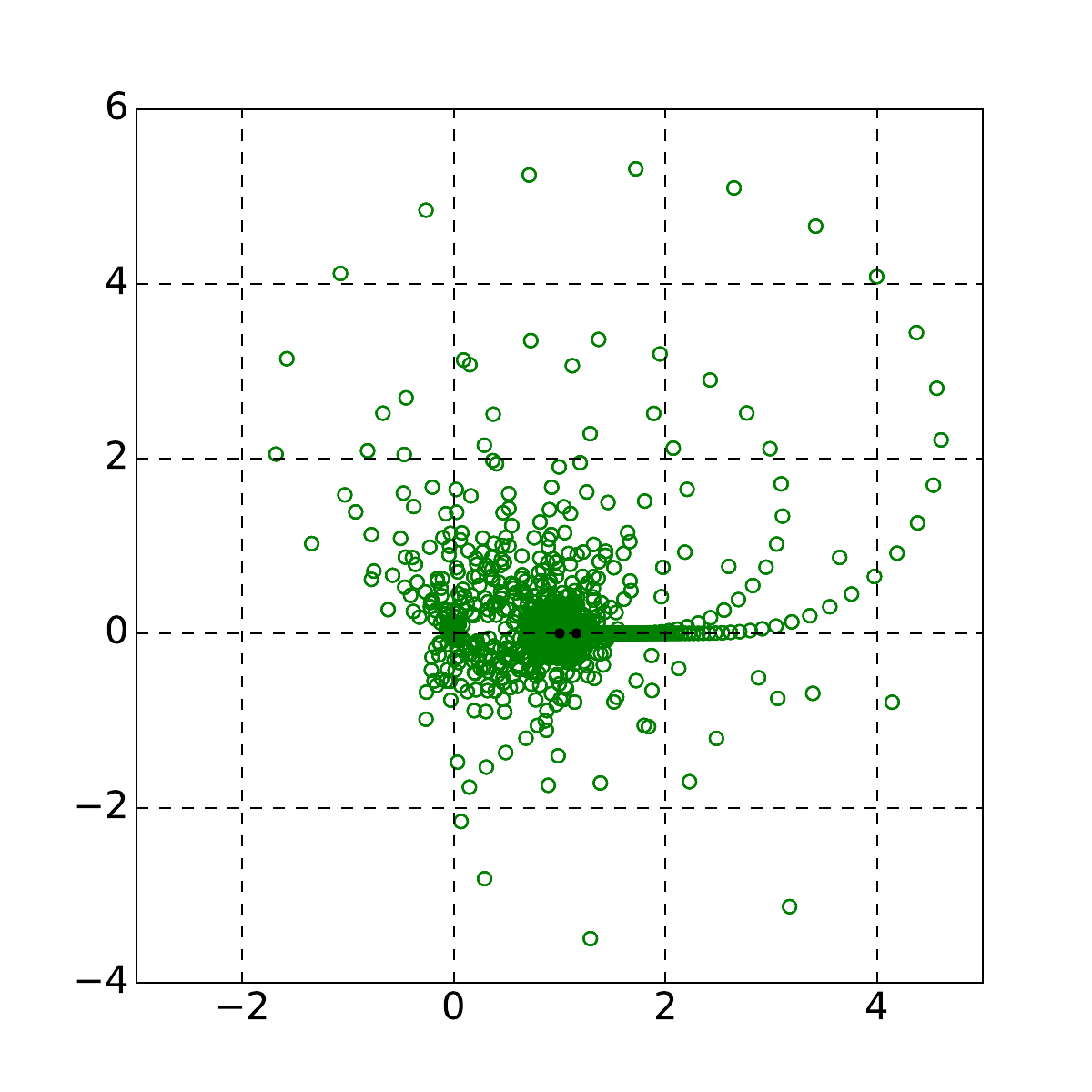}\\ \hline
 \end{tabular}
\caption{Eigenvalues distribution for the Ferrite case of the scaled versions of $\mrm{MTF}_{\loc}$ (red), $\mrm{MTF}_{\loc}^2$ (blue) and $\rev{\IB_{(\kappa,\mu)}} \cdot \mrm{MTF}_{\loc}$ (green).}
\label{tab:EigenvalueDistributionPrec}
\end{center} 
\end{table} 
\subsection{Preconditioning: Numerical Experiments}\label{Numerical}
We consider a partition of the sphere into disjoint planar triangles and we assembly the local MTF operator and the preconditioners with the open-source Galerkin boundary element library Bempp 3.3 \cite{SAB15}. A Bempp notebook server is easily accessible through Docker\footnote{https://bempp.com/download/}. The meshes and the simulations are fully reproducible as a Python Notebook\footnote{https://github.com/pescap/mtf}. Bempp allows for Calderón-based preconditioning through barycentric refinement and has the Buffa-Christiansen (BC) function basis implemented (refer to \cite{SBE17} and the references therein). We apply mass preconditioning to all formulations---strong form in Bempp---in order to obtain matrices whose condition number is bounded with the meshsize and study the eigenvalues and condition numbers of the induced operators along with the restarted GMRes(20) convergence \cite{SS86} and solver times. We set relative tolerance of GMRes(20) to $10^{-8}$ and perform simulations through a $4$ GB RAM per core, $64$ bit Linux server.

We decide to focus on the LF and HF Teflon scattering problem. For the LF (resp.~HF) case, we consider two meshes corresponding to precisions of $r_0=10$ and $r_1=15$ elements per wavelength, referred to as cases $N_0$ and $N_1$, with $N_0<N_1$ the size of the induced linear systems for $\mrm{MTF}_{\loc}$ (detailed in Table \ref{tab:simParameters}). The incident wave is a plane wave polarized along $z$-axis and traveling at a $\theta=\frac{\pi}{4}$-angle. The LF case is assembled in dense mode while we use hierarchical matrices for the HF case with ACA compression with relative tolerance $10^{-3}$. To begin with, we summarize the parameters for each case in Table \ref{tab:simParameters}. To provide an extension of the results, we introduce the HF scattering of a complex shape, namely the unit Fichera Cube---the unit cube with a reentrant corner, referred to as (HF)$_\square$.

\begin{minipage}{0.5\linewidth}
\begin{table}[H]
\renewcommand\arraystretch{1.2}
\begin{center}
\footnotesize
\begin{tabular}{|c|c|c|c|c|} \hline   
 \multicolumn{2}{|c|}{Case}&  (LF) & (HF) & (HF)$_\square$\\ \hline
 \multicolumn{2}{|c|}{$f$ (Hz)} & $50$ MHz &\multicolumn{2}{c|}{ $300$ MHz } \\ \hline\hline
 \multicolumn{2}{|c|}{$\kappa_0$} & $1.05$ &\multicolumn{2}{c|}{ $6.29$ } \\ \hline 
 \multicolumn{2}{|c|}{$\kappa_1$} & $1.52$ & \multicolumn{2}{c|}{ $9.11$ } \\ \hline 
 $r_0$ & $N_0$ & $1,380$ &  $38,856$ & $23,052$\\\cline{1-1}\hline
 $r_1$ & $N_1$ & $2,856$ & $90,276$ &$47,544$ \\\hline
\end{tabular}
\end{center}
\caption{Parameters chosen for the simulations.}
\label{tab:simParameters}
\end{table}
\end{minipage}
\begin{minipage}{0.5\linewidth}
  \hspace{1cm}
  \includegraphics[width=4cm]{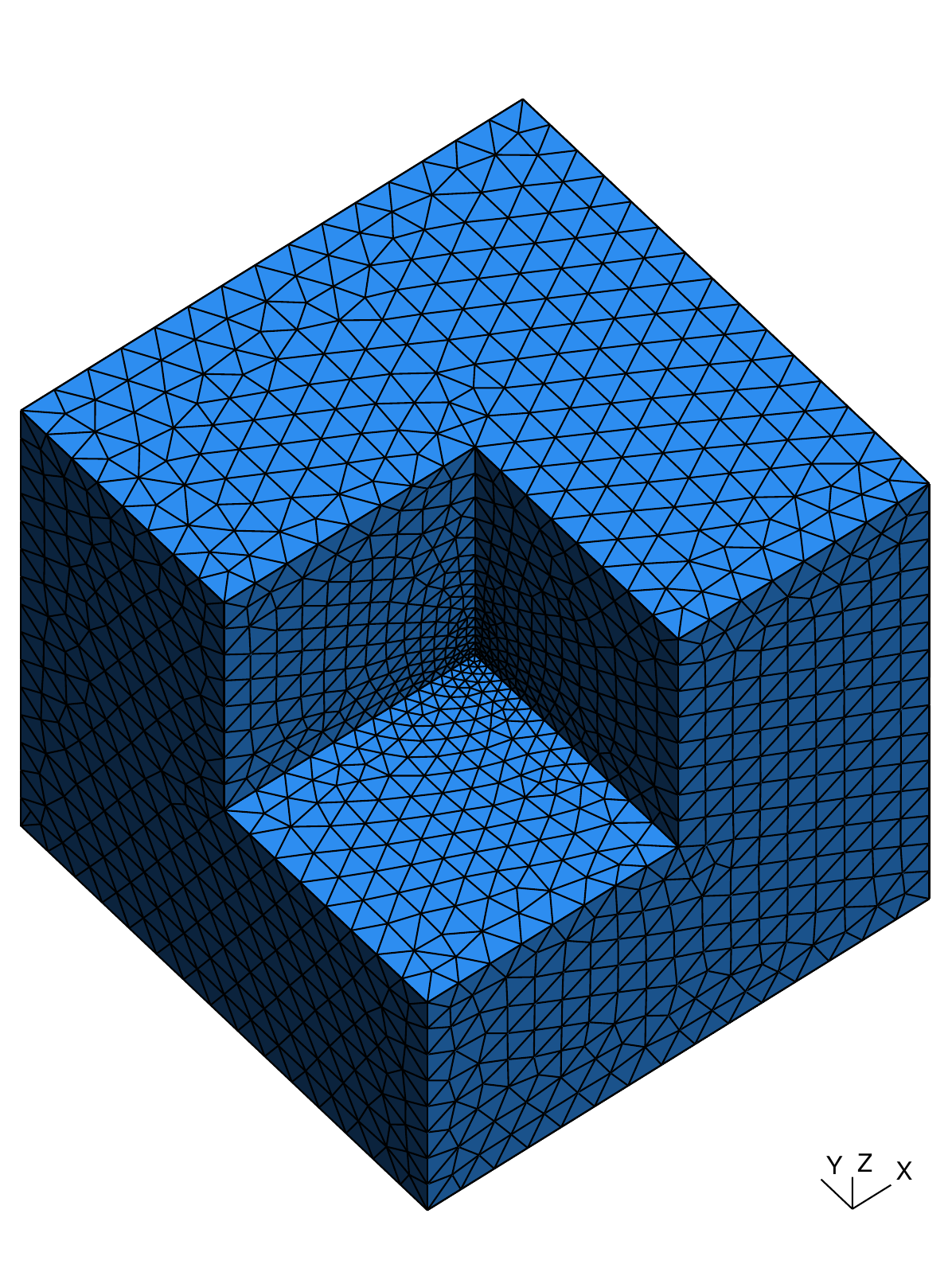}
\end{minipage}

\quad\\[5pt]
For comparison purposes, we solve the problem with the preconditioned Single-Trace-Formulation (STF) as a reference \cite{CAF11}. We verify that both the local MTF and STF lead to the same current densities. For example, the LF case with $N_0$ leads to a $1.38\%$ (resp.~$0.563\%$) relative error in $\bL^2$-norm for the exterior Dirichlet (resp.~Neumann) trace. Similarly, we obtain $0.767\%$ and $0.525\%$ for (HF)$_\square$ with $N_1$. These relative errors remain the same for all preconditioners.
\begin{table}[H]
\renewcommand\arraystretch{1.2}
\begin{center}
\footnotesize
\begin{tabular}{|c||c|c|c|c||c|c|c|c||c|c|c|c|} \hline   
Case & \multicolumn{4}{c||}{(LF)} & \multicolumn{4}{c||}{(HF)}  & \multicolumn{4}{c|}{(HF)$_\square$}\\ \hline
  Parameter& \multicolumn{2}{c|}{$n_\textup{iter}$}& \multicolumn{2}{c||}{cond.}& \multicolumn{2}{c|}{$n_\textup{iter}$}& \multicolumn{2}{c||}{$t_\textup{solve}$} & \multicolumn{2}{c|}{$n_\textup{iter}$} & \multicolumn{2}{c|}{$t_\textup{solve}$}  \\ \hline
  $N$& $N_0$ & $N_1$ & $N_0$ & $N_1$& $N_0$ & $N_1$& $N_0$ & $N_1$& $N_0$ & $N_1$& $N_0$ & $N_1$\\ \hline\hline 
 $\mrm{STF}^2$&  $8$& $8$ & $1.45\e08$  & $7.75\e07$  & $20$ &$20$& $193.2$&$612.9$ & $16$ &$16$ & $93.0$&$266$\\\hline 
 $\mrm{MTF}_{\loc}$&$28$& $27$ & $4.95\e 12$ & $5.00 \e 12$ & $108$ & $108$ & $606.8$ & $1738$ & $68$&$68$&$216$&$593$\\\hline 
 $\mrm{MTF}^2_{\loc}$&$12$&$12$ & $1.15\e 10$& $9.71 \e 09$  & $44$ & $44$& $443.1$ & $1385$ & $34$ &$34$ &$202$&$570$\\\hline 
 $\rev{\IB_{(\kappa,\mu)}} \cdot \mrm{MTF}_{\loc}$&$9$&$9$ & $5.96\e10$ & $6.07 \e 10$ & $21$ & $20$ & $1008$ & $2858$  & $17$&$17$ & $451$&$1255$\\\hline 
\end{tabular}
\end{center} 
\caption{Overview of the results for all cases. We detail $n_\textup{iter}$ (resp.~$t_\textup{solve}$) the number of iterations (resp.~total solver times in seconds) of GMRes(20) along with cond., the spectral condition number.}
\label{tab:res}
\end{table}
In Table \ref{tab:res}, we provide an extensive summary of the results. To begin with, we focus on the rows corresponding to (LF). We remark that all formulations show excellent and mesh stable convergence for GMRes(20) (see $n_\textup{iter}$). As a remark, in all cases considered in this section, the unpreconditioned GMRes(20) failed to converge. The condition number (``cond.'' row in Table \ref{tab:res}) for all formulations has a relatively high magnitude but remains stable with mesh (it even slightly betters for the three first cases with increasing $N$).

Table \ref{tab:Eigenvalue Distribution Discrete} displays the eigenvalues of the resulting matrices for (LF) and both values of $N$. We obtain a similar distribution to that expected from spectral analysis (see Fig.~\ref{fig:eigenvalues_prec}). Furthermore, the eigenvalues distribution is highly independent of the meshwidth, confirming the results presented so far. 

\begin{figure}[H]
  \centering
  \includegraphics[width=1\linewidth]{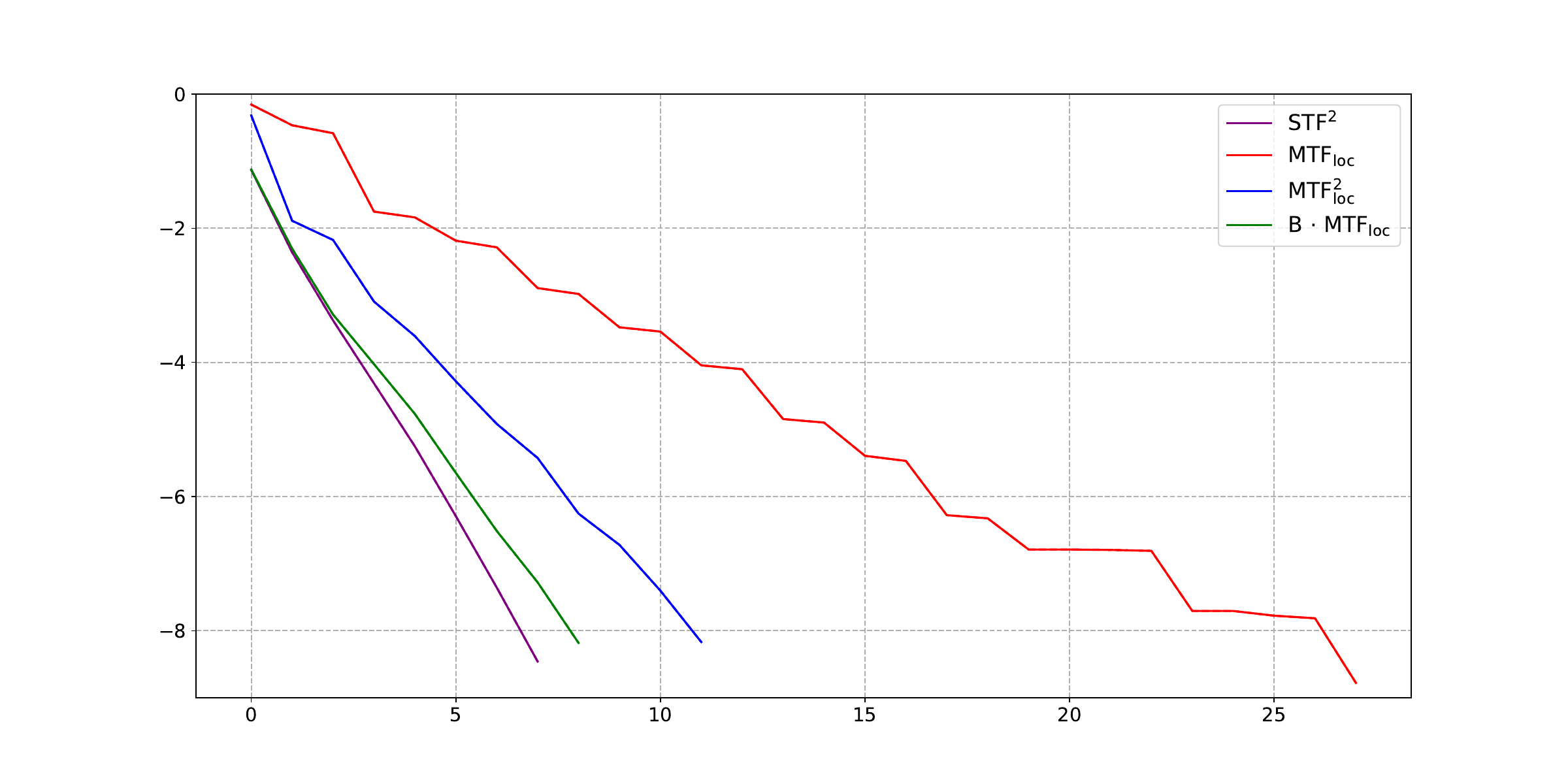}
\caption{LF: Residual convergence of GMRes(20) for $\mrm{MTF}_{\loc}$ (red), $\mrm{MTF}_{\loc}^2$ (blue), $\rev{\IB_{(\kappa,\mu)}} \cdot \mrm{MTF}_{\loc}$ (green) and $\mrm{STF}^2$ (purple) obtained with $N_0$ (dashed line) and $N_1$ (solid line).}
\label{fig:LF_eigenvalues_prec}
\end{figure}

\begin{figure}[H]
  \centering
  \includegraphics[width=1\linewidth]{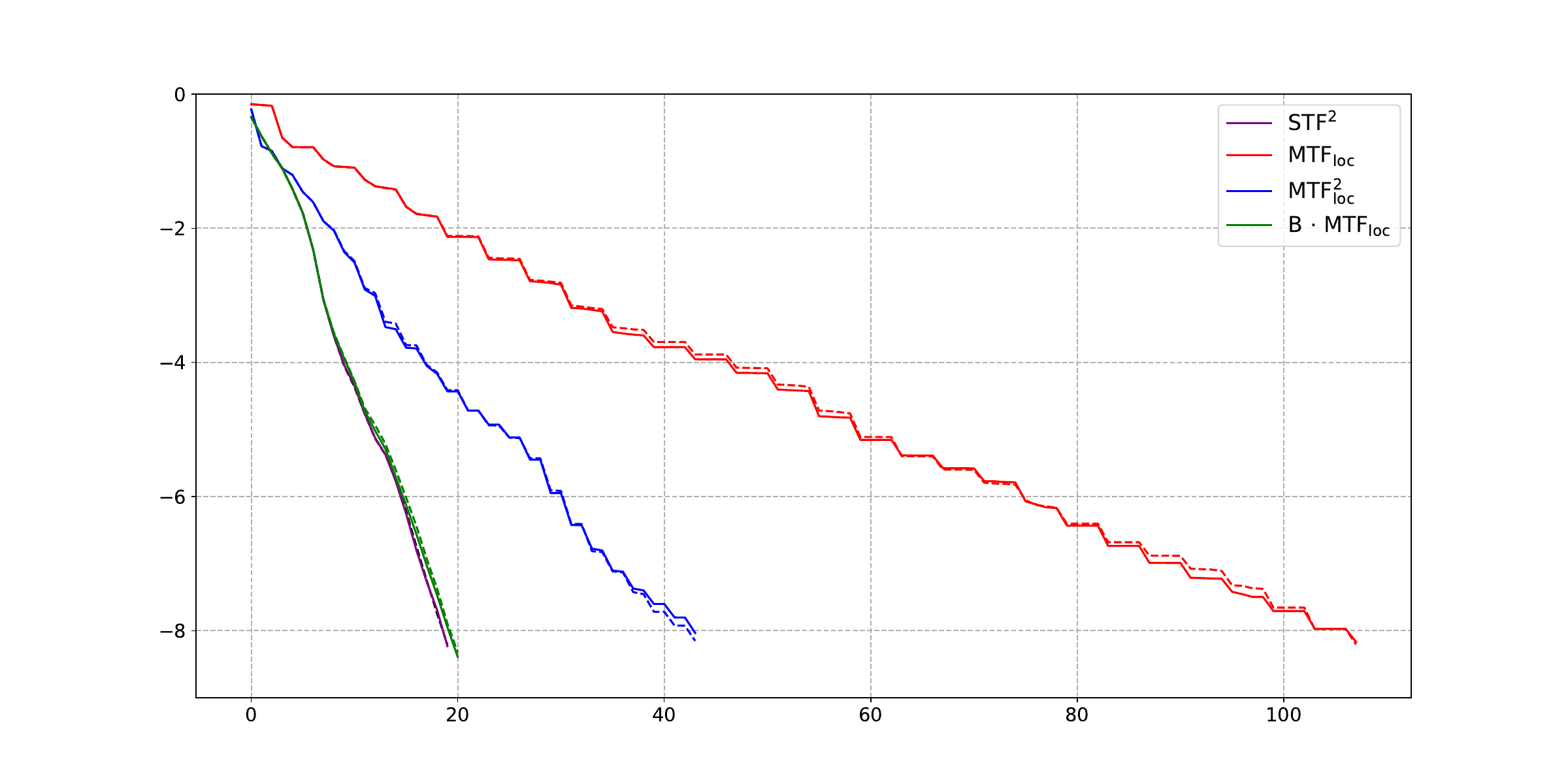}
\caption{HF: Residual convergence of GMRes(20) for $\mrm{MTF}_{\loc}$ (red), $\mrm{MTF}_{\loc}^2$ (blue), $\rev{\IB_{(\kappa,\mu)}}\cdot \mrm{MTF}_{\loc}$ (green) and $\mrm{STF}^2$ (purple) obtained with $N_0$ (dashed line) and $N_1$ (solid line).}
\label{fig:HF_eigenvalues_prec}
\end{figure}

\begin{figure}[H]
  \centering
  \includegraphics[width=1\linewidth]{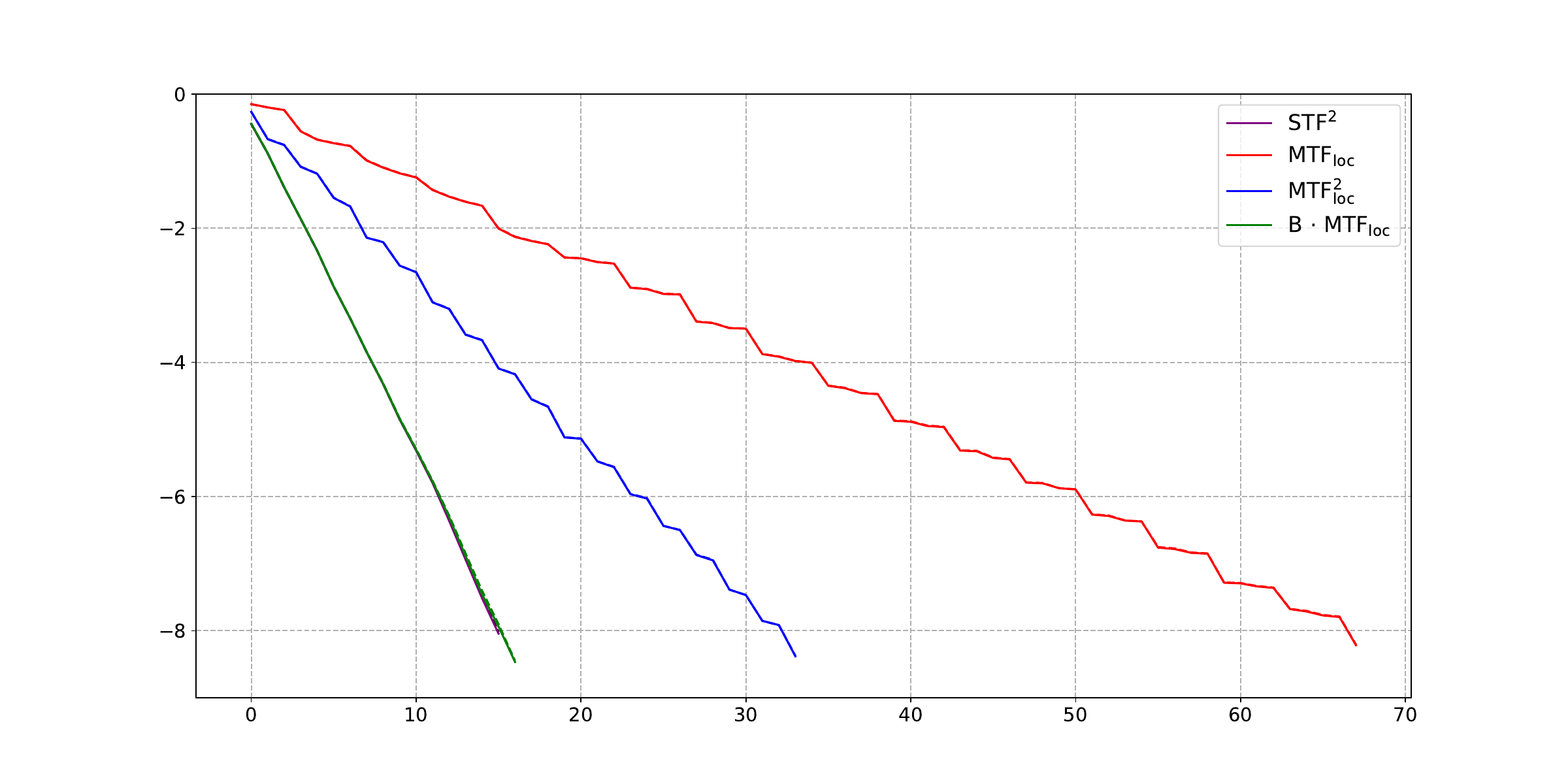}
\caption{HFS: Residual convergence of GMRes(20) for $\mrm{MTF}_{\loc}$ (red), $\mrm{MTF}_{\loc}^2$ (blue), $\rev{\IB_{(\kappa,\mu)}} \cdot \mrm{MTF}_{\loc}$ (green) and $\mrm{STF}^2$ (purple) obtained with $N_0$ (dashed line) and $N_1$ (solid line).}
\label{fig:HFS_eigenvalues_prec}
\end{figure}

\begin{table}[H]
\renewcommand\arraystretch{1.4}
\begin{center}
\footnotesize
\begin{tabular}{
    >{\centering\arraybackslash}m{1cm}
    |>{\centering\arraybackslash}m{13cm}}
\vspace{0.1cm}
 $N_0$ &  \includegraphics[width=1\linewidth]{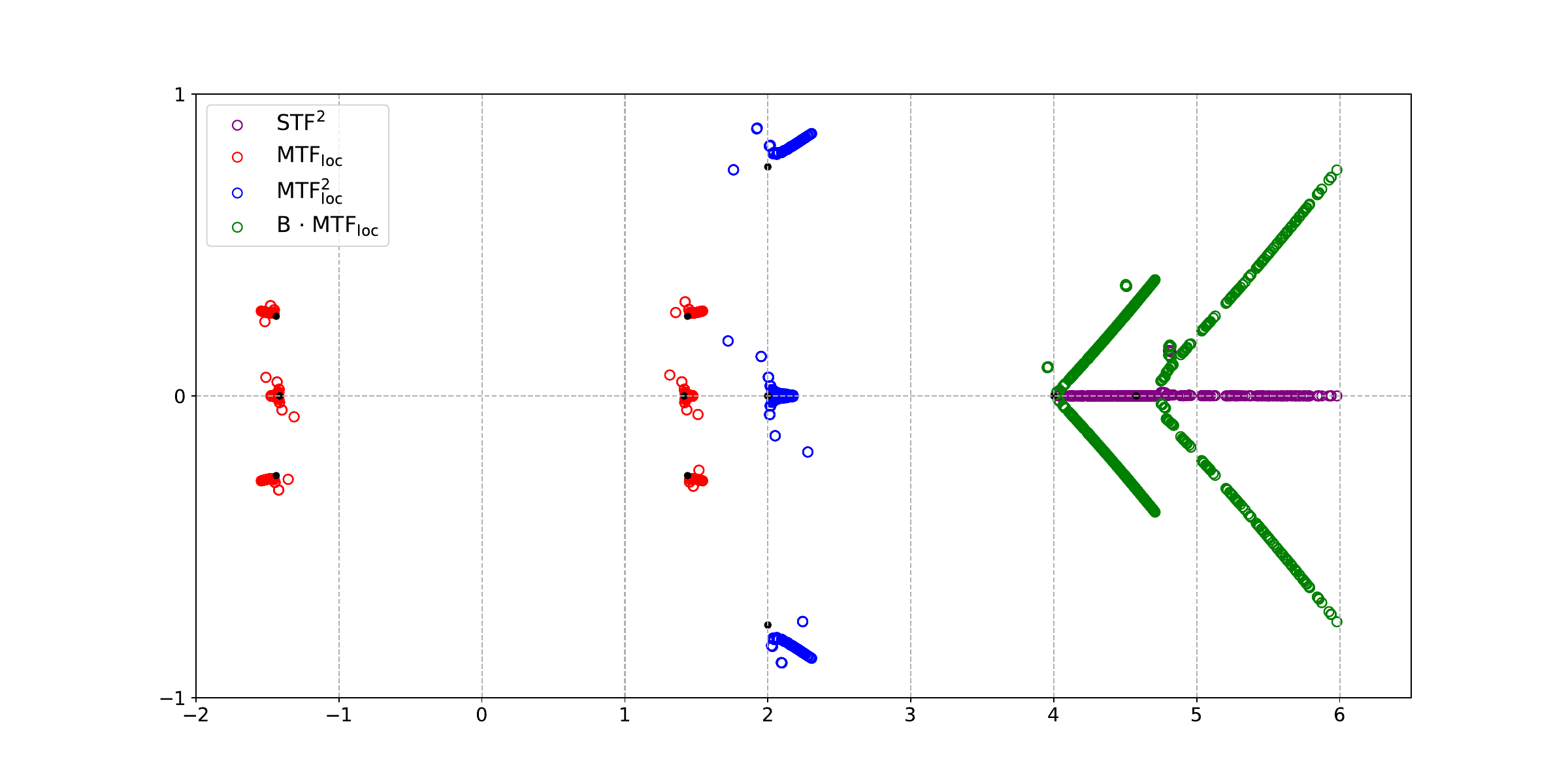}\\ \hline
  $N_1$ &  \includegraphics[width=1\linewidth]{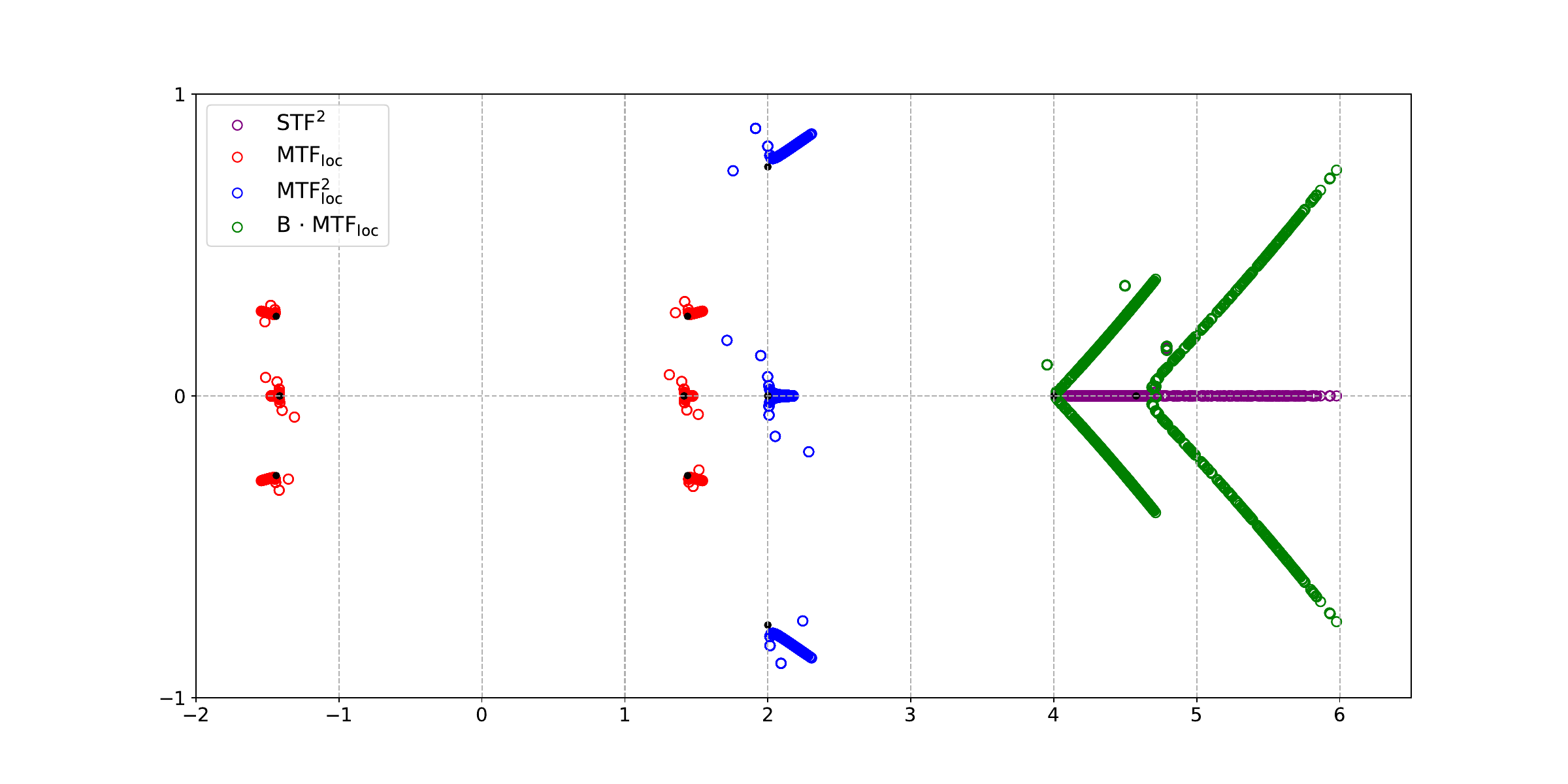} \\ \hline
 \end{tabular}
\caption{Eigenvalues distribution for the Teflon LF case of $\mrm{MTF}_{\loc}$ (red), $\mrm{MTF}_{\loc}^2$ (blue), $\rev{\IB_{(\kappa,\mu)}} \cdot \mrm{MTF}_{\loc}$ (green) and $\mrm{STF}^2$ (purple).}
\label{tab:Eigenvalue Distribution Discrete}
\end{center} 
\end{table} 
Next, in Fig.~\ref{fig:LF_eigenvalues_prec}, we plot the GMRes(20) convergence for $N_0$ (dashed line) and $N_1$ (solid line). We remark that the convergence behavior are very similar for each color. Also, the ``second-kind'' preconditioners---$\mrm{STF}^2$, $\mrm{MTF}_{\loc}^2$ and $\rev{\IB_{(\kappa,\mu)} }\cdot \mrm{MTF}_{\loc}$---outperform $\mrm{MTF}_{\loc}$. 
Afterwards, we consider both HF cases. First acknowledge that all remarks discussed before apply to those cases. In Table \ref{tab:res} (columns (HF) and (HF)$_\square$), we verify the mesh independence, as the number of iterations remains almost exactly the same with increasing mesh density. We represent the convergence of GMRes of (HF) (resp.~(HF)$_\square$) in Fig.~\ref{fig:HF_eigenvalues_prec} (resp.~Fig.~\ref{fig:HFS_eigenvalues_prec}). Again, the second-kind preconditioners outperform simple mass matrix preconditioning in terms of convergence of GMRes(20), as well as in solver time (refer to columns $t_\textup{solve}$ in Table \ref{tab:res}). Concerning mesh independence, the convergence curves for (HF)$_\square$ in Fig.~\ref{fig:HFS_eigenvalues_prec} are almost superposed, evidencing a strong mesh independence for all cases.

To finish, $\rev{\IB_{(\kappa,\mu)}} \cdot \mrm{MTF}_{\loc}$ appears to converge around twice as fast as $\mrm{MTF}_{\loc}^2$ (see $n_\textup{iter}$ in Table \ref{tab:res}) at the cost of $3$ versus $2$ matrix-vector products per iteration, benefiting a priori to $\mrm{MTF}_{\loc}^2$. More surprisingly, $\rev{\IB_{(\kappa,\mu)}} \cdot \mrm{MTF}_{\loc}$ converges at the same rate as $\mrm{STF}^2$ for the HF cases despite a two-fold increasing in degrees of freedom due to the use of multi-trace space. 

\section{Concluding remarks}\label{sec:Conclusion}
This article paves the way to show  well-posedness of the local MTF applied to electromagnetic wave scattering. \rev{As pointed out in Section \ref{sec:Introduction}, the results in this paper can be generalised to other smooth surfaces by compact perturbation.} Further research includes theoretical and numerical results for multiple domains and domains with junction or triple points. Concerning the MTF linear system preconditioning, our research hints at applying the fast preconditioning technique of \cite{EIJH19,EIJH21} to produce lower requirements of matrix-vector products for the preconditioners (applicable to both $\mrm{MTF}_{\loc}^2$ and $\rev{\IB_{(\kappa,\mu)}}\cdot\mrm{MTF}_{\loc}$). \rev{Fast convergence for the Fichera cube showed the applicability to complex (non-smooth) shapes.} Finally, the novel preconditioner $\rev{\IB_{(\kappa,\mu)}}$ proposed here paves a way toward high order inverse approximation of operators and Calder\'on-based polynomial preconditioners \cite{F91}.




\bibliography{ref-v3}





\end{document}